\documentclass[12pt]{article}

\usepackage[default]{jasa_harvard}    % for formatting citations in text
\usepackage{JASA_manu}
\usepackage{tocbibind}
\usepackage[toc,page]{appendix}

\newcommand {\ctn}{\citeasnoun} % change to \citet if using natbib

\usepackage{graphicx,subfigure,amsmath,latexsym,amssymb}
\usepackage{float,epsfig,multirow,rotating,times}
\usepackage{upgreek,wrapfig}
\usepackage{comment}

\newcommand{\bx}{\bm{x}}
\newcommand{\bX}{\boldsymbol{X}}

\newtheorem{theorem}{Theorem}

\newtheorem{corollary}[theorem]{Corollary}

\newtheorem{definition}[theorem]{Definition}

\newcommand{\e}{\ensuremath{\epsilon}}

\newcommand{\be}{\boldsymbol\e}
\newcommand{\bm}{\mathbf}

\numberwithin{equation}{section}
\numberwithin{algo}{section}
\numberwithin{table}{section}
\numberwithin{figure}{section}

\usepackage{setspace}
\usepackage[mathscr]{euscript}
\usepackage[margin=1in]{geometry}
\begin{document}
\normalsize

\title{\vspace{-0.8in}
A Short Note on Almost Sure Convergence of Bayes Factors in the General Set-Up}
\author{Debashis Chatterjee, Trisha Maitra and Sourabh Bhattacharya\thanks{
Debashis Chatterjee and Trisha Maitra are PhD students and Sourabh Bhattacharya 
is an Associate Professor in
Interdisciplinary Statistical Research Unit, Indian Statistical
Institute, 203, B. T. Road, Kolkata 700108.
Corresponding e-mail: sourabh@isical.ac.in.}}
\date{\vspace{-0.5in}}
\maketitle%

\begin{abstract}
%In this article we derive the almost sure convergence theory of Bayes factor in the general set-up that includes even
%dependent data and misspecified models, as a simple application of a result of \ctn{Shalizi09} to a well-known 
%identity satisfied by the Bayes factor.

Although there is a significant literature on the asymptotic theory of Bayes factor, the set-ups considered
are usually specialized and often involves independent and identically distributed data. Even in such specialized
cases, mostly weak consistency results are available. 
In this article, for the first time ever, we derive the almost sure convergence theory of Bayes factor in the general set-up that includes even
dependent data and misspecified models. 
Somewhat surprisingly, the key to the proof of such a general theory is a simple application of a result of \ctn{Shalizi09} to a well-known
identity satisfied by the Bayes factor.
\\[2mm]
{\it {\bf Keywords:} Bayes factor convergence; Kullback-Leibler divergence; Posterior consistency.}
 
\end{abstract}

\section{Introduction}
\label{sec:intro}

Bayes factors are well-established in the Bayesian literature for the purpose of model comparison. Briefly,
given data $\bX_n=\{X_1,X_2,\ldots,X_n\}$, where $n$ is the sample size, consider the problem of comparing 
any two models $\mathcal M_1$ and $\mathcal M_2$
associated with parameter spaces $\Theta_1$ and $\Theta_2$, respectively. For $i=1,2$, let the likelihoods, priors 
and the marginal densities for the two models be $L_n(\theta_i|\mathcal M_i)=f_{\theta_i}(\bX_n|\mathcal M_i)$, 
$\pi(\theta_i|\mathcal M_i)$ and $m(\bX_n|\mathcal M_i)=\int_{\Theta_i} L_n(\theta_i|\mathcal M_i)\pi(d\theta_i|\mathcal M_i)$,
respectively. Then the Bayes factor of model $\mathcal M_1$ against $\mathcal M_2$ is given by
\begin{equation}
B^{(12)}_n=\frac{m(\bX_n|\mathcal M_1)}{m(\bX_n|\mathcal M_2)}.
\label{eq:bf}
\end{equation}
Thus, $B^{(12)}_n$ can be interpreted as the quantification of the evidence of model $\mathcal M_1$
against model $\mathcal M_2$, given data $\bX_n$. A comprehensive account of Bayes factors is provided
in \ctn{Kass95}. 

The asymptotic study of Bayes factor involves investigation of limiting properties of $B^{(12)}_n$ as $n$ goes to infinity.
In particular, it is essential to guarantee the consistency property that $B^{(12)}_n$ goes to infinity almost surely 
when $\mathcal M_1$ is the better model and to zero almost surely when $\mathcal M_2$ is the better model. 
It is also important to obtain the rate of convergence of the Bayes factor. In the case of
independent and identically distributed ($iid$) data, a relevant result is provided in \ctn{Walker04a}
and \ctn{Walker04b}. Such strong ``almost sure" convergence results are rare however, even when the data are 
independent but not identically distributed. Recently, \ctn{Maitra16a} obtained a strong, general result 
when the data are independent but not identically distributed and applied it to time-varying covariate and drift function selection in the context of systems 
of stochastic differential equations (see also \ctn{Maitra16b} for further application of Bayes factor
asymptotics in stochastic differential equations). The other existing works on Bayes factor asymptotics 
are problem specific and even in such particular set-ups strong consistency results are seldom available
(but see, for example, \ctn{Dawid92}, \ctn{Kundu14}, \ctn{Choi15}). For a comprehensive review of Bayes 
factor consistency, see \ctn{Chib16}. 

We are interested in more general frameworks where the data may be dependent and where the possible models
are perhaps all misspecified. We are not aware of any existing work on Bayes factor asymptotics in this direction.  
However, posterior convergence has been addressed by \ctn{Shalizi09}, and indeed,
Theorem 2 of \ctn{Shalizi09} combined with a well-known identity satisfied by Bayes factors, holds the key to 
an elegant almost sure convergence result for the Bayes factor. The result depends explicitly on the average
Kullback-Leibler divergence between the competing and the true models, even in such a general set-up.  
Here it is important to emphasize that although \ctn{Chib16} is essentially a review paper, the authors demonstrate
for the first time with a specific example of nested models that the identity satisfied by the Bayes factor may be 
exploited to prove weak consistency of the latter, and provide general discussion regarding ``in probability" 
Bayes factor convergence assuming that the identity is satisfied by the Bayes factor.

The rest of this article is structured as follows. In Section \ref{sec:general_setup}, based on \ctn{Shalizi09} we describe
the general setting for our Bayes factor investigation, and provide the result of \ctn{Shalizi09} on which
our main result on Bayes factor hinges. In Section \ref{sec:bf} we provide our results on Bayes factor convergence.
We make concluding remarks in Section \ref{sec:conclusion}.
Additional details are provided in the online supplement, whose sections have the prefix ``S-" when referred to in this paper. 

\section{The general set-up for model comparison using Bayes factors}
\label{sec:general_setup}
Following \ctn{Shalizi09}, let us consider a probability space $(\Omega,\mathcal F,P)$,  sequence
of random variables $\{X_1,X_2,\ldots\}$ taking values in the measurable space $(\aleph,\mathcal X)$,
having infinite-dimensional distribution $P$. %Denoting $\left\{X_t;t=1,\ldots,T\right\}$ by $\bX_T$,
In other words, the distribution $P$ is an infinite-dimensional distribution since it is the joint distribution of 
infinitely many random variables corresponding to a valid stochastic process. 
As guaranteed by Kolmogorov's consistency result (see, for example, \ctn{Billingsley95}, \ctn{Schervish95}), 
all finite-dimensional distributions associated with $P$ can be obtained by 
marginalizing over the remaining (infinite number of) variables.
The theoretical development requires no restrictive assumptions on $P$ such as it being a product measure, Markovian, 
or exchangeable, thus paving the way for great generality. 

Let $\mathcal F_n=\sigma(\bX_n)$ denote the natural filtration, that is, the $\sigma$-algebra generated by $\bX_n$. 
%$\bX_T=\{X_1,\ldots,X_T\}$, $T$ being the sample size.
Also, let the distributions of the processes adapted to $\mathcal F_n$ be denoted by $F_{\theta}$, where
$\theta$ takes values in a measurable space $(\Theta,\mathcal T)$. Here $\theta$ denotes the hypothesized probability measure
associated with the unknown distribution of $\{X_1,X_2,\ldots\}$ and $\Theta$ is the set of hypothesized
probability measures. 
In other words, assuming that $\theta$ is the
infinite-dimensional distribution of the stochastic process $\{X_1,X_2,\ldots\}$, $F_{\theta}$ denotes the 
$n$-dimensional marginal distribution associated with $\theta$; $n$ is suppressed for ease of notation.
For parametric models, the probability 
measure $\theta$ corresponds to a probability density with respect to some dominating measure (such as Lebesgue
or counting measure) and consists of finite number of parameters. 
For nonparametric models, $\theta$ is usually associated with an infinite number of parameters and may not 
have a density with respect to $\sigma$-finite measures.

As in \ctn{Shalizi09}, we assume that $P$ and all the $F_{\theta}$ are dominated by a common measure
with densities $p$ and $f_\theta$, respectively. In \ctn{Shalizi09} and in our
case, the assumption that $P\in\Theta$ is not required so that all possible models are allowed to be misspecified.
%Indeed, \ctn{Shalizi09} provides an example of such misspecification where the true model $P$ is not Markov 
%but all the hypothesized models indexed by $\theta$ are $k$-th order stationary binary Markov models, for $k=1,2,\ldots$. 
%As shown in \ctn{Shalizi09}, the results of posterior convergence hold even in the case of such misspecification, 
%essentially because the true model can be approximated by the $k$-th order Markov models belonging to $\Theta$.

Given a prior $\pi$ on $\theta$, we assume that the posterior distributions $\pi(\cdot|\bX_n)$ are dominated by a common
measure for all $n\geq 1$; abusing notation, we denote the density at $\theta$ by $\pi(\theta|\bX_n)$.

%\begin{itemize}
%\item[(A7)] The sets $\mathcal G_T$ of (A5) and (A6) can be chosen such that for any set $A$ with $\pi(A)>0$, 
%\begin{equation}
%h\left(\mathcal G_T\cap A\right)\rightarrow h\left(A\right),
%\label{eq:A7}
%\end{equation}
%as $T\rightarrow\infty$.
%\end{itemize}
Let $L_n(\theta)=f_{\theta}(\bX_n)$ be the likelihood and $p_n=p(\bX_n)$ be the marginal density of $\bX_n$ under the true model $P$. 
Then following the notation of \ctn{Shalizi09}, for $A\subseteq\Theta$, let
%\begin{equation}
%h(\theta)=\underset{T\rightarrow\infty}{\lim}~\frac{1}{n}E\left[\log\left\{\frac{p_n}{L_n(\theta)}\right\}\right].
%\label{eq:A2}
%\end{equation}
\begin{align}
h(\theta)&=\underset{T\rightarrow\infty}{\lim}~\frac{1}{n}E\left[\log\left\{\frac{p_n}{L_n(\theta)}\right\}\right];\label{eq:A2}\\
h\left(A\right)&=\underset{\theta\in A}{\mbox{ess~inf}}~h(\theta);\label{eq:h2}\\
J(\theta)&=h(\theta)-h(\Theta);\label{eq:J}\\
J(A)&=\underset{\theta\in A}{\mbox{ess~inf}}~J(\theta),\label{eq:J2}
\end{align}
where, for any function $g:\Theta\mapsto\mathbb R$, where $\mathbb R$ is the real line, 
$$\underset{\theta\in A}{\mbox{ess~inf}}~g(\theta)=
\sup\left\{r\in\mathbb R:g(\theta)>r,~\mbox{for almost all}~\theta\in A\right\},$$ 
is the essential infimum of $g$ over the set $A$. Here ``almost all" is with respect to the prior distribution. 
In words, essential infimum is the greatest lower bound which holds with prior probability one. 
With the above notations, six assumptions are used to prove Theorem 2 of \ctn{Shalizi09}. We provide the six assumptions
in Section S-1 of the supplement, which we refer to as (A1)--(A6). Below we furnish Theorem 2 of \ctn{Shalizi09} which shall
play the key role for our purpose of deriving almost sure convergence of Bayes factors.

\begin{theorem}[Theorem 2 of \ctn{Shalizi09}]
\label{theorem:shalizi}
Consider assumptions (A1)--(A6). Then for all $\theta$ such that $\pi(\theta)>0$,
\begin{equation}
\underset{n\rightarrow\infty}{\lim}~\frac{1}{n}\log\left[\pi(\theta|\bX_n)\right]=-J(\theta),
\label{eq:shalizi}
\end{equation}
almost surely with respect to the true model $P$, where $J(\theta)$ is given by (\ref{eq:J}). 
\end{theorem}
%We shall use the above theorem to derive almost sure convergence of Bayes factors.

\section{Convergence of Bayes factors}
\label{sec:bf}

For the model comparison problem using Bayes factors, we now assume the likelihoods and the priors of all the competing models
satisfy (A1)--(A6), in addition to satisfying that $P$ and all the $F_{\theta}$ for $\theta\in\Theta_1\cup\Theta_2$ 
have densities with respect to a common $\sigma$-finite measure.
We also assume that for $i=1,2$, the posterior $\pi(\cdot|\bX_n,\mathcal M_i)$ associated with model $\mathcal M_i$ 
is dominated by the prior $\pi(\cdot|\mathcal M_i)$, which is again absolutely continuous with respect to
some appropriate $\sigma$-finite measure. These latter assumptions ensure that up to the normalizing constant, 
the posterior density %$\pi(\cdot|\bX_T,\mathcal M_i)$ 
associated with $\mathcal M_i$ is factorizable into 
the prior density times the likelihood.
Indeed, for any $\theta_i\in\Theta_i$, 
\begin{equation}
\log\left[m(\bX_n|\mathcal M_i)\right]=\log\left[L_n(\theta_i|\mathcal M_i)\right]+\log\left[\pi(\theta_i|\mathcal M_i)\right]
-\log\left[\pi(\theta_i|\bX_n,\mathcal M_i)\right].
\label{eq:logmarginal}
\end{equation}
%where $\pi(\theta_i|\bX_T,\mathcal M_i)$ denotes the posterior density of $\theta_i$ given $\bX_T$ and model $\mathcal M_i$.
Hence, the logarithm of the Bayes factor is given, for any $\theta_1\in\Theta_1$ and $\theta_2\in\Theta_2$, 
by (see, for example, \ctn{Chib95}, \ctn{Chib16})
\begin{align}
\log\left[B^{(12)}_n\right]&=\log\left[\frac{L_n(\theta_1|\mathcal M_1)}{L_n(\theta_2|\mathcal M_2)}\right]
+\log\left[\frac{\pi(\theta_1|\mathcal M_1)}{\pi(\theta_2|\mathcal M_2)}\right]
-\log\left[\frac{\pi(\theta_1|\bX_n,\mathcal M_1)}{\pi(\theta_2|\bX_n,\mathcal M_2)}\right],\notag
%&=\log\left[\frac{R_T(\theta|\mathcal M_1)}{R_T(\theta|\mathcal M_2)}\right]
%+\log\left[\frac{\pi(\theta|\mathcal M_1)}{\pi(\theta|\mathcal M_2)}\right]
%-\log\left[\frac{\pi(\theta|\bX_T,\mathcal M_1)}{\pi(\theta|\bX_T,\mathcal M_2)}\right],
%\label{eq:logbf}
\end{align}
so that
\begin{align}
\frac{1}{n}\log\left[B^{(12)}_n\right]&=\frac{1}{n}\log\left[R_n(\theta_1|\mathcal M_1)\right]
-\frac{1}{n}\log\left[R_n(\theta_2|\mathcal M_2)\right]\notag\\
&\qquad +\frac{1}{n}\log\left[\pi(\theta_1|\mathcal M_1)\right]
-\frac{1}{n}\log\left[\pi(\theta_2|\mathcal M_2)\right]\notag\\
&\qquad-\frac{1}{n}\log\left[\pi(\theta_1|\bX_n,\mathcal M_1)\right]
+\frac{1}{n}\log\left[\pi(\theta_2|\bX_n,\mathcal M_2)\right],
\label{eq:logbf2}
\end{align}
where, for $i=1,2$, $R_n(\theta_i|\mathcal M_i)=\frac{L_n(\theta_i|\mathcal M_i)}{p_n}$.

Now let $J_i(\theta_i)=h_i(\theta_i)-h_i(\Theta_i)$, where $h_i(\theta_i)$ is defined as in (\ref{eq:A2}) with 
$L_n(\theta)$ replaced with $L_n(\theta_i|\mathcal M_i)$, and 
$h_i\left(A\right)=\underset{\theta_i\in A_i}{\mbox{ess~inf}}~h_i(\theta_i)$, for any $A_i\subseteq\Theta_i$.
Assumption (A3) then yields 
\begin{equation}
\underset{n\rightarrow\infty}{\lim}~\frac{1}{n}\log \left[R_n(\theta_i|\mathcal M_i)\right]=-h_i(\theta_i),
\label{eq:R_conv}
\end{equation}
almost surely, and assuming that both the models and their associated priors satisfy assumptions (A1)--(A6), it follows using 
Theorem \ref{theorem:shalizi} that for $i=1,2$, 
\begin{equation}
\underset{n\rightarrow\infty}{\lim}~\frac{1}{n}\log\left[\pi(\theta_i|\bX_n,\mathcal M_i)\right]=-J_i(\theta_i),
\label{eq:post_conv}
\end{equation}
almost surely.

Assuming that for $i=1,2$, $\pi(\theta_i|\mathcal M_i)>0$ for all $\theta_i\in\Theta_i$, 
note that $\frac{1}{n}\log\left[\pi(\theta_i|\mathcal M_i)\right]\rightarrow 0$ as 
$n\rightarrow\infty$, so that it follows using (\ref{eq:logbf2}), (\ref{eq:R_conv}) and (\ref{eq:post_conv}), 
%and using the convention $h_i(\theta_i)-\left[h_i(\theta_i)-h_i(\Theta_i)\right]=h_i(\Theta_i)$ when $h_i(\theta_i)=\infty$, 
that
\begin{equation}
\underset{n\rightarrow\infty}{\lim}~\frac{1}{n}\log \left[B^{(12)}_n\right] = -\left[h_1(\Theta_1)-h_2(\Theta_2)\right],
\label{eq:bf_conv}
\end{equation}
almost surely with respect to $P$.
We formalize this main result in the form of the following theorem:
\begin{theorem}[Bayes factor convergence]
\label{theorem:bf_convergence}
Assume that for $i=1,2$, the competing models $\mathcal M_i$ satisfy assumptions (A1)--(A6), with parameter
spaces $\Theta_i$, in addition to satisfying that $P$ and all the $F_{\theta}$ for $\theta\in\Theta_1\cup\Theta_2$ 
have densities with respect to a common $\sigma$-finite measure.
We also assume that the posterior associated with $\mathcal M_i$ is dominated by the prior, 
which is again absolutely continuous with respect to some appropriate $\sigma$-finite measure, and  
that the priors satisfy $\pi(\theta_i|\mathcal M_i)>0$ for all $\theta_i\in\Theta_i$.
Then (\ref{eq:bf_conv}) holds almost surely with respect to the true infinite-dimensional probability measure $P$.
\end{theorem}
Since assumption (A3) is used directly for convergence of the likelihood ratios, it is perhaps desirable to 
consider sufficient conditions that ensure (A3). Such sufficient conditions, as noted in \ctn{Shalizi09}, can be found
in \ctn{Algoet88} and \ctn{Gray90}. Necessary and sufficient conditions for (A3) to hold has more recently been
established in (\ctn{Harrison08}). However, in our experience, (A3) is usually easy to verify; see Section S-2 of the
supplement; see also \ctn{Maitra16c}.

Theorem \ref{theorem:bf_convergence} provides an elegant convergence result for Bayes factors, explicitly in terms of
differences between average Kullback-Leibler divergences between the competing and the true models. That such a result
holds in the general set-up that includes even dependent data and misspecified models, is very encouraging. 
%even though it follows trivially
%from Theorem 2 of Shalizi applied directly to the well-known identity of the Bayes factor. 
Indeed, we are not aware of any such result in the general set-up, although in the $iid$ situation \ctn{Walker04a}
and \ctn{Walker04b} prove strong convergence of Bayes factor in terms of Kullback-Leibler divergences, taking misspecification
into account. Theorem \ref{theorem:bf_convergence} readily leads to the following corollaries.

\begin{corollary}[Consistency of Bayes factor]
\label{corollary:bf_consistency}
Without loss of generality, let $\mathcal M_1$ be the correct model and $\mathcal M_2$ be incorrect. Then 
$L_n(\theta_1|\mathcal M_1)=p_n$ for all $\theta_1\in\Theta_1$, so that $h_1(\theta_1)=0$ for
all $\theta_1\in\Theta_1$, implying that $h_1(\Theta_1)=0$. On the other hand, $h_2(\Theta_2)>0$, so that by 
Theorem \ref{theorem:bf_convergence},
$\underset{n\rightarrow\infty}{\lim}~\frac{1}{n}\log \left[B^{(12)}_n\right]=h_2(\Theta_2)$. In other words, $B^{(12)}_n\rightarrow\infty$ 
exponentially fast, confirming consistency of the Bayes factor. If $\mathcal M_1$ is not necessarily the correct model
but is a better model than $\mathcal M_2$ in the sense that its average Kullback-Leibler divergence $h_1(\Theta_1)$
is smaller than $h_2(\Theta_2)$, then again $B^{(12)}_n\rightarrow\infty$ exponentially fast, guaranteeing consistency.
\end{corollary}

\begin{corollary}[Selection among a countable class of models]
\label{corollary:model_selection}
Theorem \ref{theorem:bf_convergence} and Corollary \ref{corollary:bf_consistency} make it explicit that if the class of competing 
models is countable and contains the true model, it is selected by the Bayes factor, otherwise Bayes factor 
selects the model for which the average Kullback-Leibler divergence from the true model is minimized among the (countable) 
class of misspecified models, provided that the infimum is attained by some model. 
%However, for an infinite class of models $\left\{\mathcal M_s:s\in\mathcal S\right\}$,
%where $\mathcal S$ is an infinite set, $\underset{\mathcal S}{\inf}~h_s(\Theta_s)$ may be zero, even if the true
%model is not included in the considered class of models, and the infimum may not be attained by any model.
\end{corollary}

\begin{corollary}[The case when two or more models are asymptotically correct]
\label{corollary:both correct}
For simplicity let us consider two models $\mathcal M_1$ and $\mathcal M_2$ as before with parameter
spaces $\Theta_1$ and $\Theta_2$ respectively.
From Theorem \ref{theorem:bf_convergence} it follows that $\frac{1}{n}\log\left[ B^{(12)}_n\right]\rightarrow 0$
almost surely if and only if $h_1(\Theta_1)=h_2(\Theta_2)$, that is, if and only if both the models
are asymptotically correct in the average Kullback-Leibler sense. 
Note that the zero limit of $\frac{1}{n}\log\left[ B^{(12)}_n\right]$ is the only logical limit here 
since any non-zero limit would lead the Bayes factor
to lend infinitely more support to one model compared to the other even though both the competing models are
correct asymptotically.
The situation of zero limit of $\frac{1}{n}\log\left[ B^{(12)}_n\right]$ may arise
in the case of comparisons between nested models or when testing parametric versus nonparametric models.
In these cases even though both the competing models are correct asymptotically, one may be a much larger model.
For reasons of parsimony it then makes sense to choose the model with smaller dimensionality. If both the models
are infinite-dimensional, for example, when comparing two sets of basis functions, then model combination seems 
to be the right step.
\end{corollary}

In Section S-2 of the supplement we illustrate Theorem \ref{theorem:bf_convergence} with an example with
autoregressive models of the first order ($AR(1)$ models) comparing (asymptotically) stationary versus nonstationary models
when the true model is (asymptotically) stationary. We show that asymptotically the Bayes factor heavily favours the (asymptotically) stationary model.

In Corollary \ref{corollary:both correct}, we have referred to comparisons with nonparametric models. 
However, recall that the results of Shalizi require the true model $P$ 
and all the postulated models $F_{\theta}$ to have densities with respect to a common dominating measure, and also the posteriors
$\pi(\cdot|\bX_n)$ to be dominated by a common reference measure for all $n\geq 1$. These conditions are typically satisfied
by parametric models, but not necessarily by nonparametric models. Indeed, in the case of the traditional nonparametric
Bayesian analysis using the Dirichlet process prior, there is no parametric form of the likelihood as there is no density
of the data $\bX_n$ under this nonparametric set-up. Also, the prior is not dominated by any $\sigma$-finite measure, 
and so does not have any density. In other words, not all nonparametric models lead to posteriors that can be factorized 
as proportional to prior times likelihood, as our Bayes factor treatment requires.  
However, as we clarify in Section S-3 of the supplement with a series of various examples of 
nonparametric Bayesian set-ups, in general
the aforementioned factorization of the posterior holds in Bayesian nonparametrics and the domination requirements
of Shalizi also hold in general. However, we emphasize that we did not yet verify assumptions (A1)--(A6) for these
cases, as we reserve this task for our future paper to be communicated elsewhere.

\section{Conclusion}
\label{sec:conclusion}
In this article, we have obtained an elegant almost sure convergence result for Bayes factors in the general set-up
where the data may be dependent and where all possible models are allowed to be misspecified. To our knowledge,
this is a first-time effort in this direction. Interestingly, in spite of the importance of the result, it 
follows rather trivially from Shalizi's result on posterior consistency applied to the identity satisfied by Bayes factors.
We assert that although similar results can be shown to hold in simpler set-ups 
(see \ctn{Walker04a} and \ctn{Walker04b} for the $iid$ set-up and \ctn{Maitra16a} for the independent
and non-identical set-up) and perhaps under specific models, our contribution
is a proof of a strong convergence result under a very general set-up that has not been considered before. 

The generality of our result will enable Bayes factor based asymptotic comparisons of various models in various set-ups, 
for example, $k$-th order Markov models, hidden Markov models, spatial Markov random
field models, models based on dependent systems of stochastic differential equations, parametric versus nonparametric
models in the dependent data setting (\ctn{Ghosal08} consider the $iid$ set-up and study ``in-probability" 
convergence of Bayes factor comparing specific finite and infinite-dimensional models).
dependent versus independent model set-ups, to name only a few. Moreover, even in the $iid$
data contexts, the existing Bayes factor asymptotic results for the specific problems are usually not directly based
on Kullback-Leibler divergence. Since our result directly make use of Kullback-Leibler divergence in any set-up, it is much more
appealing from this perspective compared to the existing results. 

In our future endeavors, we shall explore the effectiveness of our result in various specific set-ups, along with comparisons
with existing results whenever applicable.

\section*{Acknowledgment}
We are sincerely grateful to the Editor-in-Chief, the Associate Editor and the referee whose detailed constructive comments have led
to significant improvement of our manuscript.

\newpage

\renewcommand\thefigure{S-\arabic{figure}}
\renewcommand\thetable{S-\arabic{table}}
\renewcommand\thesection{S-\arabic{section}}

\setcounter{section}{0}
\setcounter{figure}{0}
\setcounter{table}{0}

\begin{center}
{\bf \Large Supplementary Material}
\end{center}

\section{Assumptions of Shalizi in the context of posterior consistency}
\label{subsec:assumptions_shalizi}

\begin{itemize}
\item[(A1)] 
%Letting $L_n(\theta)=f_\theta(\bX_n)$ be the likelihood and
%$p_n=p(\bX_n)$ be the marginal distribution of $\bX_n$ under the true model $P$, 
Assume that the following likelihood ratio
\begin{equation}
R_n(\theta)=\frac{L_n(\theta)}{p_n}
\label{eq:R_T}
\end{equation}
%Assume that $R_n(\theta)$ 
is $\mathcal F_n\times \mathcal T$-measurable for all $n\geq 1$.
\end{itemize}

\begin{itemize}
\item[(A2)] For every $\theta\in\Theta$, the Kullback-Leibler divergence rate $h(\theta)$
%\begin{equation}
%h(\theta)=\underset{T\rightarrow\infty}{\lim}~\frac{1}{n}E\left[\log\left\{\frac{p_n}{L_n(\theta)}\right\}\right].
%\label{eq:A2}
%\end{equation}
exists (possibly being infinite) and is $\mathcal T$-measurable.
Note that in the $iid$ set-up, $h(\theta)$ reduces to the Kullback-Leibler divergence between the true and the
hypothesized model, so that %(\ref{eq:A2}) 
$h(\theta)$ may be regarded as a generalized Kullback-Leibler divergence measure.
\end{itemize}

\begin{itemize}
\item[(A3)] For each $\theta\in\Theta$, the generalized or relative asymptotic equipartition property holds, and so,
almost surely with respect to $P$,
\begin{equation}
\underset{n\rightarrow\infty}{\lim}~\frac{1}{n}\log \left[R_n(\theta)\right]=-h(\theta).
\label{eq:equipartition}
\end{equation}
%where $h(\theta)$ is given by (\ref{eq:A2}).

Roughly, the terminology 
``asymptotic equipartition" refers to dividing up $\log \left[R_n(\theta)\right]$ into $n$ factors for large $n$ 
such that all the factors are asymptotically equal. Again, considering the $iid$ scenario helps clarify this point,
as in this case each factor converges to the same Kullback-Leibler divergence between the true and the postulated model. 
With this understanding note that the purpose of condition (A3) is to ensure that relative to the true distribution, 
the likelihood of each $\theta$ decreases to zero exponentially
fast, with rate being the Kullback-Leibler divergence rate (\ref{eq:equipartition}). 
\end{itemize}

\begin{itemize}
\item[(A4)] 
Let $I=\left\{\theta:h(\theta)=\infty\right\}$. 
The prior $\pi$ on $\theta$ satisfies $\pi(I)<1$.
Failure of this assumption entails extreme misspecification of almost all the hypothesized models $f_{\theta}$ relative
to the true model $p$. With such extreme misspecification, posterior consistency is not expected to hold; see
\ctn{Shalizi09} for details.
\end{itemize}

%Following the notation of \ctn{Shalizi09}, for $A\subseteq\Theta$, let
%\begin{align}
%h\left(A\right)&=\underset{\theta\in A}{\mbox{ess~inf}}~h(\theta);\label{eq:h2}\\
%J(\theta)&=h(\theta)-h(\Theta);\label{eq:J}\\
%J(A)&=\underset{\theta\in A}{\mbox{ess~inf}}~J(\theta),\label{eq:J2}
%\end{align}
%where, for any function $g:\Theta\mapsto\mathbb R$, where $\mathbb R$ is the real line, 
%$$\underset{\theta\in A}{\mbox{ess~inf}}~g(\theta)=
%\sup\left\{r\in\mathbb R:g(\theta)>r,~\mbox{for almost all}~\theta\in A\right\},$$ 
%is the essential infimum of $g$ over the set $A$. Here ``almost all" is with respect to the prior distribution. 
%In words, essential infimum is the greatest lower bound which holds with prior probability one. 
\begin{itemize}
\item[(A5)] There exists a sequence of sets $\mathcal G_n\rightarrow\Theta$ as $n\rightarrow\infty$ 
such that: %along with $\pi(\mathcal G_T)>0$
\begin{enumerate}
\item[(1)] $h\left(\mathcal G_n\right)\rightarrow h\left(\Theta\right)$, as $n\rightarrow\infty$.
\item[(2)]
\begin{equation}
\pi\left(\mathcal G_n\right)\geq 1-\alpha\exp\left(-\beta n\right),~\mbox{for some}~\alpha>0,~\beta>2h(\Theta);
\label{eq:A5_1}
\end{equation}
\item[(3)]The convergence in (A3) is uniform in $\theta$ over $\mathcal G_n\setminus I$.
\end{enumerate}

The sets $\mathcal G_n$ can be loosely interpreted as the sieves corresponding to the method of sieves (\ctn{Geman82})
such that the behaviour of the likelihood ratio and the posterior on the sets $\mathcal G_n$ 
essentially carries over to $\Theta$. 
This can be anticipated from the first and the second parts of the assumption; the second part ensuring in particular
that the parts of $\Theta$ on which the log likelihood ratio may be ill-behaved have exponentially small prior probabilities.
The third part is more of a technical condition that is useful in proving posterior convergence through the sets
$\mathcal G_n$. For further details, see \ctn{Shalizi09}.
\end{itemize}
For each measurable $A\subseteq\Theta$, for every $\delta>0$, there exists a random natural number $\tau(A,\delta)$
such that
\begin{equation}
n^{-1}\log\left[\int_{A}R_n(\theta)\pi(\theta)d\theta\right]
\leq \delta+\underset{n\rightarrow\infty}{\lim\sup}~n^{-1}
\log\left[\int_{A}R_n(\theta)\pi(\theta)d\theta\right],
\label{eq:limsup_2}
\end{equation}
for all $n>\tau(A,\delta)$, provided 
$\underset{n\rightarrow\infty}{\lim\sup}~n^{-1}\log\left[\int_{A}R_n(\theta)\pi(\theta)d\theta\right]<\infty$.
%$\mathbb I_A$ denotes the indicator function of the set $A$.
Regarding this, the following assumption has been made by Shalizi:
\begin{itemize}
\item[(A6)] The sets $\mathcal G_n$ of (A5) can be chosen such that for every $\delta>0$, the inequality
$n>\tau(\mathcal G_n,\delta)$ holds almost surely for all sufficiently large $n$.

To understand the essence of this assumption, note that for almost every data set $\{X_1,X_2,\ldots\}$
there exists $\tau(\mathcal G_n,\delta)$ such that (\ref{eq:limsup_2}) holds with $A$ replaced by $\mathcal G_n$
for all $n>\tau(\mathcal G_n,\delta)$. Since $\mathcal G_n$ are sets with large enough prior probabilities,
the assumption formalizes our expectation that $R_n(\theta)$ decays fast enough on $\mathcal G_n$.
%so that $\tau(\mathcal G_n,\delta)$ is nearly stable in the sense that it is not only finite but also not 
%significantly different for different data sets when $n$ is large. 
See \ctn{Shalizi09} for more detailed explanation.

\end{itemize}

\section{Illustration of our result on Bayes factor with competing $AR(1)$ models}
\label{sec:bf_illustration}

Let the true model $P$ stand for the following $AR(1)$ model: 
\begin{equation}
x_t=\rho_0 x_{t-1}+\epsilon_t,~t=1,2,\ldots,
\label{eq:true_ar1}
\end{equation}
where $x_0\equiv 0$, $|\rho_0|<1$ and $\epsilon_t\stackrel{iid}{\sim}N(0,\sigma^2_0)$, for $t=1,2,\ldots$.
We assume the competing models $\mathcal M_1$ and $\mathcal M_2$ to be of the same form as
(\ref{eq:true_ar1}) but with the true parameter $\rho_0$ replaced with the unknown parameters $\rho_1$ and $\rho_2$, respectively,
such that $|\rho_1|<1$ and $\rho_2\in (-1,1)^c\cap\mathbb S$, where $(-1,1)^c$ denotes complement of $(-1,1)$ and $\mathbb S$ is some compact set containing $[-1,1]$. 
For model $\mathcal M_i$; $i=1,2$, we assume that 
$x_0\equiv 0$ and $\epsilon_t\stackrel{iid}{\sim}N(0,\sigma^2_i)$; $t=1,2,\ldots$. 
For simplicity of illustration we assume for the time being that $\sigma_1$ and $\sigma_2$ are known, that is, 
$\sigma_1=\sigma_2=\sigma_0$, 
but see Section \ref{subsec:sigma_unknown} where we allow $\sigma_1$ and $\sigma_2$ to be unknown.
Thus, we are interested in comparing (asymptotically) stationary and nonstationary $AR(1)$ models where the true $AR(1)$ model is (asymptotically) stationary.  
Note that $\Theta_1=(-1,1)$ and $\Theta_2=(-1,1)^c\cap\mathbb S$. We consider priors $\pi(\cdot|\mathcal M_i)$; $i=1,2$, both of which
have densities with respect to the Lebesgue measure.
Let us first verify assumptions (A1)--(A6) with respect to $\mathcal M_1$. All the probabilities and expectations
below are with respect to the true model $P$. Notationally, in this time series context we denote the sample size by the more natural notation $T$ rather than $n$. 

\subsection{Verification of (A1) for $\mathcal M_1$}
\label{subsec:A1}
Note that
\begin{equation}
\log R_T(\rho_1)=
%\left(\frac{\rho_0-\rho_1}{\sigma^2}\right)
%\left[\sum_{t=1}^Tx_tx_{t-1}-\left(\sum_{t=1}^Tx^2_{t-1}\right) \left(\rho_0+\rho_1\right)\right].
\left(\frac{\rho_0-\rho_1}{\sigma^2_0}\right)
\left[\left(\sum_{t=1}^Tx^2_{t-1}\right) \left(\frac{\rho_0+\rho_1}{2}\right)-\sum_{t=1}^Tx_tx_{t-1}\right].
\label{eq:BF_R_T}
\end{equation}
Thanks to continuity it is clear that $R_T(\rho_1)$ is %clearly a random function of $\bX_T$ and $\rho_1$ and hence is 
$\mathcal F_T\times\mathcal T$ measurable. In other words, (A1) holds.

\subsection{Verification of (A2) for $\mathcal M_1$}
\label{subsec:A2}
It is easy to verify that under the true model $P$ the autocovariance function is given by
\begin{equation}
Cov(x_{t+h},x_t)\sim\frac{\sigma^2_0\rho^h_0}{1-\rho^2_0};~h\geq 0,
\label{eq:cov_true}
\end{equation}
where for any two sequences $\{a_t\}_{t=1}^{\infty}$ and $\{b_t\}_{t=1}^{\infty}$, $a_t\sim b_t$ stands for $a_t/b_t\rightarrow 1$ as $t\rightarrow\infty$.
%see, for example, \ctn{Shumway11},
This leads to
\begin{align}
E\left[\log R_T(\rho_1)\right]&=-\left(\frac{\rho_1-\rho_0}{\sigma^2_0}\right)
\left[\left(\sum_{t=1}^TE\left(x^2_{t-1}\right)\right) \left(\frac{\rho_1+\rho_0}{2}\right)
-\sum_{t=1}^TE\left(x_tx_{t-1}\right)\right]
\notag\\
&\sim-\left(\rho_1-\rho_0\right)
\left[\frac{(T-1)\left(\rho_1+\rho_0\right)}{2(1-\rho^2_0)}-\frac{(T-1)\rho_0}{(1-\rho^2_0)}\right],\notag
\end{align}
so that
\begin{equation*}
\frac{E\left[\log R_T(\rho_1)\right]}{T}\rightarrow -\frac{\left(\rho_1-\rho_0\right)^2}{2(1-\rho^2_0)},
~\mbox{as}~T\rightarrow\infty.
\end{equation*}
In other words, (A2) holds, with 
\begin{equation}
h_1(\rho_1)=\frac{\left(\rho_1-\rho_0\right)^2}{2(1-\rho^2_0)}.
\label{eq:h_1}
\end{equation}

\subsection{Verification of (A3) for $\mathcal M_1$}
\label{subsec:A3}
Rather than proving pointwise almost sure convergence of $\frac{\log R_T(\rho_1)}{T}$ to $-h_1(\rho_1)$,
we prove the stronger result of almost sure uniform convergence in our example. 
Indeed, note that
\begin{align}
&\underset{|\rho_1|<1}{\sup}~\left|\frac{\log R_T(\rho_1)}{T}+h_1(\rho_1)\right|\notag\\
&=\underset{|\rho_1|<1}{\sup}~\left|\frac{\rho_1-\rho_0}{\sigma^2_0}\right|\times
\left|\left(\frac{\sum_{t=1}^Tx^2_{t-1}}{T}\right) \left(\frac{\rho_1+\rho_0}{2}\right)-\frac{\sum_{t=1}^Tx_tx_{t-1}}{T}-
\frac{\sigma^2_0\left(\rho_1-\rho_0\right)}{2(1-\rho^2_0)}\right|\notag\\
&\leq\underset{|\rho_1|\leq 1}{\sup}~\left|\frac{\rho_1-\rho_0}{\sigma^2_0}\right|\times
\left|\left(\frac{\sum_{t=1}^Tx^2_{t-1}}{T}\right) \left(\frac{\rho_1+\rho_0}{2}\right)-\frac{\sum_{t=1}^Tx_tx_{t-1}}{T}-
\frac{\sigma^2_0\left(\rho_1-\rho_0\right)}{2(1-\rho^2_0)}\right|\notag\\
&=\left|\frac{\hat\rho_1-\rho_0}{\sigma^2_0}\right|\times
\left|\left(\frac{\sum_{t=1}^Tx^2_{t-1}}{T}\right) \left(\frac{\hat\rho_1+\rho_0}{2}\right)-\frac{\sum_{t=1}^Tx_tx_{t-1}}{T}-
\frac{\sigma^2_0\left(\hat\rho_1-\rho_0\right)}{2(1-\rho^2_0)}\right|
\label{eq:A3_1}\\
&\leq \kappa\left|\left(\frac{\sum_{t=1}^Tx^2_{t-1}}{T}\right) \left(\frac{\hat\rho_1+\rho_0}{2}\right)
-\frac{\sum_{t=1}^Tx_tx_{t-1}}{T}-
\frac{\sigma^2_0\left(\hat\rho_1-\rho_0\right)}{2(1-\rho^2_0)}\right|,
\label{eq:A3_2}
\end{align}
where step (\ref{eq:A3_1}) follows due to compactness of $[-1,1]$; here $\hat\rho_1\in [-1,1]$ depends upon the data.
In (\ref{eq:A3_2}), $\kappa$ is a finite positive constant greater than the bounded positive quantity 
$\left|\frac{\hat\rho_1-\rho_0}{\sigma^2_0}\right|$.  

Now observe that under $P$, the Markov chain $\left\{x_t:t=1,2,\ldots,\right\}$ is not only an asymptotically stationary process but is
also irreducible and aperiodic (for definitions, see, for example, \ctn{Meyn93} and \ctn{Robert04}). The latter two properties are easy to 
see because the chain can travel from any value in the real line to any set with positive Lebesgue measure in just one step
with positive probability. Thus, the ergodic theorem holds, so that as $T\rightarrow\infty$,
\begin{equation}
\frac{\sum_{t=1}^Tx^2_{t-1}}{T} \rightarrow %E(x^2_1)=
\frac{\sigma^2_0}{1-\rho^2_0},
\label{eq:ergodic1}
\end{equation}
almost surely with respect to $P$. To deal with $\frac{\sum_{t=1}^Tx_tx_{t-1}}{T}$, note that under $P$,
\begin{equation}
x_tx_{t-1}=\rho_0x^2_{t-1}+\epsilon_tx_{t-1},
\label{eq:decomp}
\end{equation}
and that $\left\{\epsilon_tx_{t-1}:t=2,3,\ldots\right\}$ is also
an asymptotically stationary, irreducible and aperiodic Markov chain. Hence, applying ergodic theorem to the latter Markov chain, 
we obtain, using independence of $\epsilon_t$ and $x_{t-1}$ for all $t\geq 2$,
%of $\epsilon_2$ and $x_1$,
\begin{equation}
\frac{\sum_{t=1}^T\epsilon_tx_{t-1}}{T} \rightarrow %E(\epsilon_2x_1)=E(\epsilon_2)E(x_1)=
0,
\label{eq:ergodic2}
\end{equation}
as $T\rightarrow\infty$, almost surely with respect to $P$.
It follows by combining (\ref{eq:ergodic1}), (\ref{eq:decomp}) and (\ref{eq:ergodic2}) that
\begin{equation}
\frac{\sum_{t=1}^Tx_tx_{t-1}}{T} \rightarrow \frac{\sigma^2_0\rho_0}{1-\rho^2_0},
\label{eq:ergodic3}
\end{equation}
 as $T\rightarrow\infty$, almost surely with respect to $P$.
Applying (\ref{eq:ergodic1}) and (\ref{eq:ergodic3}) to (\ref{eq:A3_2}) yields
\begin{align}
&\left|\left(\frac{\sum_{t=1}^Tx^2_{t-1}}{T}\right) \left(\frac{\hat\rho_1+\rho_0}{2}\right)-\frac{\sum_{t=1}^Tx_tx_{t-1}}{T}-
\frac{\sigma^2_0\left(\hat\rho_1-\rho_0\right)}{2(1-\rho^2_0)}\right|\notag\\
&=\left|\left(\frac{\sum_{t=1}^Tx^2_{t-1}}{T}-\frac{\sigma^2_0}{1-\rho^2_0}\right)\left(\frac{\hat\rho_1+\rho_0}{2}\right)
-\left(\frac{\sum_{t=1}^Tx_tx_{t-1}}{T}-\frac{\sigma^2_0\rho_0}{1-\rho^2_0}\right)\right|\notag\\
&\leq \left|\left(\frac{\hat\rho_1+\rho_0}{2}\right)\right|
\times\left|\frac{\sum_{t=1}^Tx^2_{t-1}}{T}-\frac{\sigma^2_0}{1-\rho^2_0}\right|
+\left|\frac{\sum_{t=1}^Tx_tx_{t-1}}{T}-\frac{\sigma^2_0\rho_0}{1-\rho^2_0}\right|\notag\\
&\rightarrow 0,
\end{align}
as $T\rightarrow\infty$, almost surely with respect to $P$.
In other words, (A3) holds and the convergence is uniform.

\subsection{Verification of (A4) for $\mathcal M_1$}
\label{subsec:A4}

In our example, (A4) holds trivially since $h_1(\rho_1)=\frac{(\rho_1-\rho_0)^2}{2(1-\rho^2_0)}$, and $|\rho|<1$
almost surely. Specifically, $\pi(I|\mathcal M_1)=0$. 

\subsection{Verification of (A5) for $\mathcal M_1$}
\label{subsec:A5}

First note that 
$h_1\left(\Theta_1\right)=\underset{\rho_1\in \Theta_1}{\mbox{ess~inf}}~h_1(\rho_1)
=\underset{\rho_1\in \Theta_1}{\mbox{ess~inf}}~\frac{(\rho_1-\rho_0)^2}{2(1-\rho^2_0)}=0$.
%Next, let $\mathcal G_T=\left\{\rho_1\in\Theta_1:|\rho_1|\leq\exp(\beta T)\right\}$, where
%$\beta>h_1(\Theta_1)=0$.
Next, let $\mathcal G_T=\Theta_1$, for $T>0$. Then (A5) (1) and (A5) (2) hold trivially.
%Clearly, $\mathcal G_T\rightarrow\Theta_1$ and $h_1(\mathcal G_T)\rightarrow 0=h_1(\Theta_1)$,
%as $T\rightarrow\infty$, so that (A5) (1) holds.
%By Markov's inequality, $\pi(\mathcal G_T)>1-\left(E|\rho_1|\right)\exp\left(-\beta T\right)$.
%That is, (A5) (2) is satisfied. 
Validation of (A5) (3) is exactly the same as our proof of uniform convergence of $\frac{\log R_T(\cdot)}{T}$
to $h_1(\cdot)$, provided in Section \ref{subsec:A3}.
Hence, (A5) is satisfied.

\subsection{Verification of (A6) for $\mathcal M_1$}
\label{subsec:A6}

Under (A1) -- (A3), which we have already verified,
it holds that (see equation (18) of \ctn{Shalizi09}) for any fixed $\mathcal G$ of the sequence $\mathcal G_T$, 
for any $\epsilon>0$ and
for sufficiently large $T$,
\begin{equation}
\frac{1}{T}\log\int_{\mathcal G}R_T(\rho_1)\pi(\rho_1|\mathcal M_1)d\rho_1
\leq -h_1(\mathcal G)+\epsilon+\frac{1}{T}\log\pi(\mathcal G|\mathcal M_1).
\label{eq:shalizi_18}
\end{equation}
It follows that $\tau(\mathcal G_T,\delta)$ is almost surely finite for all $T$ and $\delta$.
We now argue that for sufficiently large $T$, $\tau(\mathcal G_T,\delta)>T$ only finitely often
with probability one. 
By equation (41) of \ctn{Shalizi09},
\begin{align}
\sum_{T=1}^{\infty}P\left(\tau(\mathcal G_T,\delta)>T\right)
\leq \sum_{T=1}^{\infty}\sum_{m=T+1}^{\infty}
P\left(\frac{1}{m}\log\int_{\mathcal G_T}R_m(\rho_1)\pi(\rho_1|\mathcal M_1)d\rho_1>\delta-h_1(\mathcal G_T)\right).
\label{eq:A6_1}
\end{align}
%In the above, we can replace $\mathcal G_T$ with the compact set 
%Now let $\tilde{\mathcal G_T}=\left\{\rho_1\in[-1,1]\cap[-\exp(\beta T),\exp(\beta T)]\right\}$, and note that
%$\frac{1}{m}\log\int_{\mathcal G_T}R_m(\rho_1)\pi(\rho_1|\mathcal M_1)d\rho_1=
%\frac{1}{m}\log\int_{\tilde{\mathcal G_T}}R_m(\rho_1)\pi(\rho_1|\mathcal M_1)d\rho_1$.
%Abusing notation, we continue to denote this
%compact set by $\mathcal G_T$.
Since $\frac{1}{m}\log\int_{\mathcal G_T}R_m(\rho_1)\pi(\rho_1|\mathcal M_1)d\rho_1=
\frac{1}{m}\log\int_{|\rho_1|\leq 1}R_m(\rho_1)\pi(\rho_1|\mathcal M_1)d\rho_1$, by the mean value theorem for integrals,
\begin{equation}
\frac{1}{m}\log\int_{\mathcal G_T}R_m(\rho_1)\pi(\rho_1|\mathcal M_1)d\rho_1
=\frac{1}{m}\log \left[R_m(\hat\rho_T)\pi(\Theta_1|\mathcal M_1)\right]
=\frac{1}{m}\log\left[R_m(\hat\rho_T)\right],
\label{eq:mvt}
\end{equation}
for $\hat\rho_T\in[-1,1]$ depending upon the data.

%Next observe that if $(-1,1)\subset[-\exp(\beta T),\exp(\beta T)]$, then
%$h_1(\mathcal G_T)=h_1((-1,1))=0$, otherwise, $(-1,1)\supset[-\exp(\beta T),\exp(\beta T)]$, so that
%$\mathcal G_T=[-\exp(\beta T),\exp(\beta T)]$ is compact, and 
%%by compactness of $\tilde{\mathcal G_T}$, 
%hence, $h_1(\mathcal G_T)=h_1(\hat\rho_T)$, for
%$\hat\rho_T\in\mathcal G_T$. 
%That is, 
%\begin{equation}
%h_1(\mathcal G_T)=\max\left\{0,h_1(\hat\rho_T)\right\}. 
%\label{eq:h1_GT}
%\end{equation}
%Note that $h_1(\hat\rho_T)$ may also be zero.

%Since $h_1(\tilde\rho_T)\geq \max\left\{0,h_1(\hat\rho_T)\right\}=h_1(\mathcal G_T)$, 
%$$\frac{1}{m}\log\int_{\tilde{\mathcal G_T}}R_m(\rho_1)\pi(\rho_1|\mathcal M_1)d\rho_1>\delta-h_1(\mathcal G_T)$$ 
%implies, in conjunction with (\ref{eq:mvt}) and (\ref{eq:h1_GT}) that

Since $h_1(\mathcal G_T)=h_1\left((-1,1)\right)=0$, and $h_1(\hat\rho_T)\geq 0$, it follows from
$$\frac{1}{m}\log\int_{\mathcal G_T}R_m(\rho_1)\pi(\rho_1|\mathcal M_1)d\rho_1>\delta-h_1(\mathcal G_T)$$ 
and (\ref{eq:mvt}) that
$$
%\frac{1}{m}\log R_m(\tilde\rho_T)+h_1(\tilde\rho_T)\geq\delta-\frac{1}{m}\log\pi(\mathcal G_T|\mathcal M_1)>\delta.
\frac{1}{m}\log R_m(\hat\rho_T)+h_1(\hat\rho_T)>\delta+h_1(\hat\rho_T)>\delta.
$$
Thus, it follows from (\ref{eq:A6_1}), (\ref{eq:A3_2}) and (\ref{eq:decomp}), that
\begin{align}
&\sum_{T=1}^{\infty}P\left(\tau(\mathcal G_T,\delta)>T\right)\notag\\
&\leq \sum_{T=1}^{\infty}\sum_{m=T+1}^{\infty}
P\left(\left|\frac{1}{m}\log R_m(\hat\rho_T)+h_1(\hat\rho_T)\right|>\delta\right)\notag\\
&\leq \sum_{T=1}^{\infty}\sum_{m=T+1}^{\infty}
P\left(\left|\left(\frac{\sum_{t=1}^mx^2_{t-1}}{m}\right) \left(\frac{\hat\rho_T-\rho_0}{2}\right)
-\frac{\sum_{t=2}^m\epsilon_tx_{t-1}}{m}-
\frac{\sigma^2_0\left(\hat\rho_T-\rho_0\right)}{2(1-\rho^2_0)}\right|>\frac{\delta}{\kappa}\right)\notag\\
&\leq \sum_{T=1}^{\infty}\sum_{m=T+1}^{\infty}
P\left(\left|\left(\frac{\sum_{t=1}^mx^2_{t-1}}{m}\right) \left(\frac{\hat\rho_T-\rho_0}{2}\right)
-\frac{\sigma^2_0\left(\hat\rho_T-\rho_0\right)}{2(1-\rho^2_0)}\right|
>\frac{\delta}{2\kappa}\right)\label{eq:A6_0}\\
&\qquad+\sum_{T=1}^{\infty}\sum_{m=T+1}^{\infty}P\left(\left|\frac{\sum_{t=2}^m\epsilon_tx_{t-1}}{m}\right|>\frac{\delta}{2\kappa}\right).
\label{eq:A6_2}
\end{align}
%Due to (\ref{eq:ergodic1}), $\frac{\sum_{t=1}^mx^2_{t-1}}{m}\rightarrow \frac{\sigma^2_0}{1-\rho^2_0}$ as
%$m\rightarrow\infty$, almost surely
%with respect to $P$ and by (\ref{eq:ergodic2}), %the ergodic theorem also ensures that
%$\frac{\sum_{t=2}^m\epsilon_tx_{t-1}}{m}\rightarrow 0$ as $m\rightarrow\infty$, almost surely with respect to $P$.

We first show that (\ref{eq:A6_0}) is convergent. To simplify arguments, we first approximate $x_t=\sum_{k=1}^t\rho^{t-k}_0\epsilon_k$ by 
\begin{equation}
\tilde x_t=\sum_{k=t-t_0}^t\rho^{t-k}_0\epsilon_k
\label{eq:tilde_x_t}
\end{equation}
in the ``in probability" sense. In $\tilde x_t$, $t_0$ is such that, for any given $\varepsilon>0$, for $t>t_0$,
\begin{equation}
\max\left\{E\left|\epsilon_1\right|\times\frac{\rho^{t_0+1}_0}{1-\rho_0},\frac{\sigma^2_0\rho^{2(t_0+1)}_0}{1-\rho^2_0}\right\}<\varepsilon.
\label{eq:t_0_varepsilon}
\end{equation}
Since $\tilde x_t$ consists of only $t_0+1$ terms for any $t>t_0$, it is easier to handle compared to $x_t$, whose number of terms increases with $t$.
Importantly, $\tilde x_t$ and $\tilde x_{t+t_0+k}$ are independent, for any $k\geq 1$. This property, which is not possessed by $x_t$, will be instrumental for
making most of the terms zero associated with multinomial expansions required in our proceeding.

For the ``in probability" fact, note that
\begin{align}
E\left|x_t-\tilde x_t\right|\leq E\left|\epsilon_1\right|\sum_{k=1}^{t-t_0-1}\rho^{t-k}_0=E\left|\epsilon_1\right|\times\frac{\rho^{t_0+1}_0\left(1-\rho^{t-t_0-1}_0\right)}{1-\rho_0}
<\varepsilon,
\end{align}
and
\begin{align}
E\left|x_t-\tilde x_t\right|^2=\sigma^2_0\sum_{k=1}^{t-t_0-1}\rho^{2(t-k)}_0=\frac{\sigma^2_0\rho^{2(t_0+1)}_0\left(1-\rho^{2(t-t_0-1)}_0\right)}{1-\rho^2_0}
<\varepsilon,
\end{align}
due to (\ref{eq:t_0_varepsilon}).
Since $\varepsilon>0$ is arbitrary, it follows that 
\begin{equation}
\left|x_t-\tilde x_t\right|\stackrel{P}{\longrightarrow}0,~\mbox{as}~t\rightarrow\infty,
\label{eq:in_prob1}
\end{equation}
where $``\stackrel{P}{\longrightarrow}"$ indicates ``in probability" convergence.
Now, $\left|x^2_t-\tilde x^2_t\right|=\left|x_t+\tilde x_t\right|\times\left|x_t-\tilde x_t\right|$, where $x_t$ is an irreducible, aperiodic
Markov chain with mean zero Gaussian asymptotic stationary distribution with variance $\sigma^2_0/(1-\rho^2_0)$, and $\tilde x_t$ is also asymptotically Gaussian
with mean zero and variance $\sigma^2_0(1-\rho^{2(t_0+1)}_0)/(1-\rho^2_0)$. Hence, $\left|x_t+\tilde x_t\right|$ converges in probability to a finite random variable,
and because of (\ref{eq:in_prob1}), it follows from the above representation that
\begin{equation}
\left|x^2_t-\tilde x^2_t\right|\stackrel{P}{\longrightarrow}0,~\mbox{as}~t\rightarrow\infty.
\label{eq:in_prob2}
\end{equation}
It then follows from the representation 
\begin{equation*}
\left|\frac{\sum_{t=1}^Tx^2_t}{T}-\frac{\sum_{t=1}^T\tilde x^2_t}{T}\right|\leq\frac{\sum_{t=1}^T\left|x^2_t-\tilde x^2_t\right|}{T},
\end{equation*}
(\ref{eq:in_prob2}), and Theorem 7.15 of \ctn{Schervish95} that
\begin{equation}
\left|\frac{\sum_{t=1}^mx^2_t}{m}-\frac{\sum_{t=1}^m\tilde x^2_t}{m}\right|\stackrel{P}{\longrightarrow}0,~\mbox{as}~m\rightarrow\infty.
\label{eq:in_prob3}
\end{equation}

Now note that for any finite integer $p\geq 1$,
\begin{equation}
\underset{m\geq 1}{\sup}~E\left(\frac{\sum_{t=1}^mx^2_t}{m}-\frac{\sum_{t=1}^m\tilde x^2_t}{m}\right)^p
\leq 2^{p-1}\underset{m\geq 1}{\sup}~E\left(\frac{\sum_{t=1}^mx^2_t}{m}\right)^p+2^{p-1}\underset{m\geq 1}{\sup}~E\left(\frac{\sum_{t=1}^m\tilde x^2_t}{m}\right)^p.
\label{eq:ui1}
\end{equation}
Noting that the multinomial expansion $(a_1+a_2+\cdots+a_m)^p=\sum_{b_1+b_2+\cdots+b_m=p}\prod_{j=1}^ma^{b_j}_j$ (where $b_1,\ldots,b_m$ are non-negative
integers) consists of ${{m+p-1}\choose p}$ terms,
it follows using asymptotic stationarity of $x_t$ and $\tilde x_t$ that both the expectations on the right hand side of (\ref{eq:ui1}) are of the order $O(1)$, as $m\rightarrow\infty$. 
Also, since for any finite $m$, the
expectations are finite, it follows that the right hand side of (\ref{eq:ui1}) is finite, from which uniform integrability, and hence
\begin{equation}
E\left|\frac{\sum_{t=1}^mx^2_t}{m}-\frac{\sum_{t=1}^m\tilde x^2_t}{m}\right|^p\rightarrow 0,~\mbox{as}~m\rightarrow\infty,
\label{eq:ui2}
\end{equation}
follows for integers $p\geq 1$.
Hence, using binomial expansion, the Cauchy-Schwartz inequality and (\ref{eq:ui2}), it follows that
\begin{align}
&E\left|\frac{\sum_{t=1}^mx^2_t}{m}\right|^p-E\left|\frac{\sum_{t=1}^m\tilde x^2_t}{m}\right|^p\notag\\
&\qquad=E\left|\left(\frac{\sum_{t=1}^mx^2_t}{m}-\frac{\sum_{t=1}^m\tilde x^2_t}{m}\right)+\frac{\sum_{t=1}^m\tilde x^2_t}{m}\right|^p-E\left|\frac{\sum_{t=1}^m\tilde x^2_t}{m}\right|^p\notag\\
&\qquad\leq\sum_{k=0}^{p-1}{p\choose k}\left\{E\left(\left|\frac{\sum_{t=1}^mx^2_t}{m}-\frac{\sum_{t=1}^m\tilde x^2_t}{m}\right|^{2(p-k)}\right)\right\}^{1/2}
\times\left\{E\left(\left|\frac{\sum_{t=1}^m\tilde x^2_t}{m}\right|^{2k}\right)\right\}^{1/2},\notag
%&\qquad\rightarrow 0,~\mbox{as}~m\rightarrow\infty.
%\label{eq:ui3}
\end{align}
so that
\begin{align}
\frac{E\left|\frac{\sum_{t=1}^mx^2_t}{m}\right|^p}{E\left|\frac{\sum_{t=1}^m\tilde x^2_t}{m}\right|^p}-1
&\leq\sum_{k=0}^{p-1}{p\choose k}\left\{E\left(\left|\frac{\sum_{t=1}^mx^2_t}{m}-\frac{\sum_{t=1}^m\tilde x^2_t}{m}\right|^{2(p-k)}\right)\right\}^{1/2}
\times\left\{\frac{E\left(\left|\frac{\sum_{t=1}^m\tilde x^2_t}{m}\right|^{2k}\right)}{E\left|\frac{\sum_{t=1}^m\tilde x^2_t}{m}\right|^{2p}}\right\}^{1/2}\notag\\
&\qquad\rightarrow 0,~\mbox{as}~m\rightarrow\infty~\mbox{(due to (\ref{eq:ui2}))}.
\label{eq:ui4}
\end{align}

In other words, for $p\geq 1$,
\begin{equation}
E\left|\frac{\sum_{t=1}^mx^2_t}{m}\right|^p\sim E\left|\frac{\sum_{t=1}^m\tilde x^2_t}{m}\right|^p,~\mbox{as}~m\rightarrow\infty.
\label{eq:ui5}
\end{equation}
Hence, while applying Markov's inequality to the probability terms of the series (\ref{eq:A6_0}), we can replace the moments associated with $x_t$ with
those associated with $\tilde x_t$, for $m>T_0$, where $T_0$ is sufficienty large.
%It follows from (\ref{eq:in_prob3}) and (\ref{eq:ui3}) that

%Now observe that 
%\begin{align}
%&P\left(\left|\left(\frac{\sum_{t=1}^mx^2_{t-1}}{m}\right) \left(\frac{\hat\rho_T-\rho_0}{2}\right)
%-\frac{\sigma^2_0\left(\hat\rho_T-\rho_0\right)}{2(1-\rho^2_0)}\right|
%>\frac{\delta}{2\kappa}\right)\notag\\
%&\qquad\sim
%P\left(\left|\left(\frac{\sum_{t=1}^m\tilde x^2_{t-1}}{m}\right) \left(\frac{\hat\rho_T-\rho_0}{2}\right)
%-\frac{\sigma^2_0\left(\hat\rho_T-\rho_0\right)}{2(1-\rho^2_0)}\right|
%>\frac{\delta}{2\kappa}\right).
%\label{eq:asymp_equiv1}
%\end{align}
Now observe that
\begin{align}
&P\left(\left|\left(\frac{\sum_{t=1}^mx^2_{t-1}}{m}\right) \left(\frac{\hat\rho_T-\rho_0}{2}\right)
-\frac{\sigma^2_0\left(\hat\rho_T-\rho_0\right)}{2(1-\rho^2_0)}\right| >\frac{\delta}{2\kappa}\right)\notag\\
&\qquad\leq 
P\left(\left|\left(\frac{\sum_{t=1}^m\left[x^2_{t-1}-E\left(x^2_{t-1}\right)\right]}{m}\right) \left(\frac{\hat\rho_T-\rho_0}{2}\right)\right| 
>\frac{\delta}{4\kappa}\right)
\label{eq:prob1}\\
&\qquad\qquad+
P\left(\left|\frac{\hat\rho_T-\rho_0}{2}\right|\times\left|\frac{\sum_{t=1}^mE\left(x^2_t\right)}{m}-\frac{\sigma^2_0}{1-\rho^2_0}\right|>\frac{\delta}{4\kappa}\right).
\label{eq:prob2}
\end{align}
For $m>T_0$, where $T_0$ is sufficiently large, %and for $t_0$ sufficiently large, 
$\left|\frac{\hat\rho_T-\rho_0}{2}\right|\times\left|\frac{\sum_{t=1}^mE\left(x^2_t\right)}{m}-\frac{\sigma^2_0}{1-\rho^2_0}\right|<\frac{\delta}{4\kappa}$, so that
the (\ref{eq:prob2}) is exactly zero for $m>T_0$.
Using Markov's inequality for (\ref{eq:prob2}) where $m>T_0$ and replacing $x_t$ with $\tilde x_t$ in the right hand side of Markov's inequality using (\ref{eq:ui5}) we obtain
\begin{align}
&P\left(\left|\left(\frac{\sum_{t=1}^m\left[x^2_{t-1}-E\left(x^2_{t-1}\right)\right]}{m}\right) \left(\frac{\hat\rho_T-\rho_0}{2}\right)\right|>\frac{\delta}{4\kappa}\right)\notag\\
&\qquad<C\left(\frac{4\kappa}{\delta}\right)^5
\left(\frac{\hat\rho_T-\rho_0}{2}\right)^5
E\left(\frac{\sum_{t=1}^m\left[\tilde x^2_{t-1}- E\left(\tilde x^2_{t-1}\right)\right]}{m}\right)^5,
\label{eq:prob3}
\end{align}
where $C$ is a positive constant.
Now, $\left(\sum_{t=1}^m\left[\tilde x^2_{t-1}- E\left(\tilde x^2_{t-1}\right)\right]\right)^5$ admits the multinomial expansion of the form
$(a_1+a_2+\cdots+a_m)^5=\sum_{b_1+b_2+\cdots+b_m=5}\prod_{t=1}^ma^{b_t}_t$, where $a_t=\left[\tilde x^2_{t-1}- E\left(\tilde x^2_{t-1}\right)\right]$ and $b_1,\ldots,b_m$ are non-negative integers.
Observe that for any $t\geq 1$, $a_t$ and $a_{t+t_0+k}$ are independent for any $k\geq 1$, which enables factorization of $E\left(\prod_{t=1}^ma^{b_t}_t\right)$ into 
products of expectations of the independent terms. Since $E(a_t)=0$ for $t=2,\ldots,m$, the expected product term becomes zero whenever it consists of at least one term of the 
form $E(a_t)$, for any $t=2,\ldots,m$. 

For the sake of convenience, let $m=(s+1)(t_0+1)$, where $s~(\geq 1)$ is an integer. 
Let $A_l=\left\{a_t:t=(l-1)t_0+1,\ldots,l(t_0+1)\right\}$, for $l=1,\ldots,(s+1)$. Then $A_l$ and $A_{l+2+r}$ are independent sets for any integer $l\geq 1$ and any integer $r\geq 0$.

When at least one $b_t=1$, the following argument gives an upper bound on the number of ways $E\left(\prod_{t=1}^ma^{b_t}_t\right)$ can be non-zero. 
Consider selecting 5 sets, say, $\left\{A_l,A_{l+1},A_{l+2},A_{l+3},A_{l+4}\right\}$ from $\{A_1,\ldots,A_{s+1}\}$, for some $l\geq 1$. Let
$B_l=\left\{b_t:t=(l-1)t_0+1,\ldots,l(t_0+1)\right\}$ for $l=1,\ldots,(s+1)$, and consider setting one element of each of $B_{l+r}$; $r=1,\ldots,5$ to be 1 and the rest of the $b_t$'s to be zero.
Then the number of such cases, namely, $O\left((s+1)\right)$ (since $t_0$ is a constant), provides an upper bound on the number of possible ways $E\left(\prod_{t=1}^ma^{b_t}_t\right)$ can be non-zero
when at least one $b_t=1$.

%and any $5$ elements of the set 
%$B_l=\left\{b_t:t=(l-1)t_0+l,\ldots,l(t_0+1)\right\}$ are chosen to be 1 and the rest zeros. There are $O\left(s+1\right)=O(m)$ ways for this to occur, since $t_0$ is a constant. The expression
%$E\left(\prod_{t=1}^ma^{b_t}_t\right)$ may also be non-zero when $A_l$ and $A_{l+1}$ are chosen, any $r$ elements of $B_l$ and any $5-r$ elements of $B_{l+1}$ are chosen to be 1,
%for $r=1,2,3,4$, and the rest of the $b_t$ are set equal to zero. This can happen in $O\left((s+1)^2\right)=O\left(m^2\right)$ ways.
%In the same situations, $b_t$ taking the values 2 and 3 may also be envisaged. The total number of ways remains $O\left(m^2\right)$. 

Further cases of non-zero $E\left(\prod_{t=1}^ma^{b_t}_t\right)$ can occur when one of the $b_t$'s is $5$ and the rest are zeros, and when one
of the $b_t$ is $3$, another is $2$, and the rest are zeros, so that there are $m+m(m-1)=m^2$ cases with respect to such choices. 

Hence, in all there are $O\left(m^2\right)$ possible cases when 
$E\left(\prod_{t=1}^ma^{b_t}_t\right)$ is non-zero, and in the remaining cases
$E\left(\prod_{t=1}^ma^{b_t}_t\right)=0$. In other words,  
\begin{equation}
\left(\frac{4\kappa}{\delta}\right)^5\left(\frac{\hat\rho_T-\rho_0}{2}\right)^5
E\left(\frac{\sum_{t=1}^m\left[\tilde x^2_{t-1}- E\left(\tilde x^2_{t-1}\right)\right]}{m}\right)^5=O\left(m^{-3}\right),
\label{eq:order1}
\end{equation}
since $\hat\rho_T\in[-1,1]$.

Now, (\ref{eq:A6_0}) converges if and only if 
\begin{align}
%&\underset{N\rightarrow\infty}{\lim}~i
&\sum_{T=T_0}^{\infty}\sum_{m=T+1}^{\infty}
P\left(\left|\left(\frac{\sum_{t=1}^mx^2_{t-1}}{m}\right) \left(\frac{\hat\rho_T-\rho_0}{2}\right)
-\frac{\sigma^2_0\left(\hat\rho_T-\rho_0\right)}{2(1-\rho^2_0)}\right|
>\frac{\delta}{\kappa}\right)\label{eq:series1}\\
&\qquad<\infty,
\notag%\label{eq:series2}
\end{align}
for sufficiently large $T_0$. Due to %(\ref{eq:asymp_equiv1}), 
(\ref{eq:prob1}), (\ref{eq:prob2}) (which is exactly zero for $m>T_0$), (\ref{eq:prob3}) and (\ref{eq:order1}),
we see that (\ref{eq:series1}) is dominated by some 
finite positive constant times the series
\begin{align}
%E\left|\epsilon_2\right|^3\times E\left|\epsilon_1\right|^3
%\left(\frac{\kappa}{\delta}\right)^3 %\underset{N\rightarrow\infty}{\lim}~
\sum_{T=T_0}^{\infty}\sum_{m=T+1}^{\infty}\frac{1}{m^3}
&=\frac{1}{(T_0+1)^3}+\frac{1}{(T_0+2)^3}+\frac{1}{(T_0+3)^3}+\cdots\notag\\
&\qquad\qquad\qquad+\frac{1}{(T_0+2)^3}+\frac{1}{(T_0+3)^3}+\cdots\notag\\
&\qquad\qquad\qquad\qquad\qquad\qquad+\frac{1}{(T_0+3)^3}+\cdots\notag\\
& \qquad\qquad\qquad\qquad\qquad\qquad\qquad\qquad\qquad+\cdots\notag\\
&\qquad\qquad\qquad\qquad\qquad\qquad\qquad\qquad\qquad\vdots\notag\\
&=\sum_{k=1}^{\infty}\frac{k}{(T_0+k)^3}.
\label{eq:series3}
\end{align}
The series (\ref{eq:series3}) is convergent since it is bounded above by 
$\sum_{k=1}^{\infty}\frac{(T_0+k)}{(T_0+k)^3}\leq\sum_{k=1}^{\infty}\frac{1}{k^2}<\infty$.

Similar (and simpler) arguments and using the result 
\begin{equation*}
\left|\frac{\sum_{t=1}^m\epsilon_tx_{t-1}}{m}-\frac{\sum_{t=1}^m\epsilon_t\tilde x_{t-1}}{m}\right|
\leq\frac{\sum_{t=1}^m\left|\epsilon_t\right|\left|x_{t-1}-\tilde x_{t-1}\right|}{m}
\stackrel{P}{\longrightarrow}0,~\mbox{as}~m\rightarrow\infty,
%\label{eq:in_prob4}
\end{equation*}
shows that the series (\ref{eq:A6_2}) also converges.
Hence, (A6) stands verified.

Thus, (A1)--(A6) holds for $\mathcal M_1$. 

\subsection{Verification of Shalizi's conditions for model $\mathcal M_2$}
\label{subsec:verify_M2}
We now verify the same set of conditions for $\mathcal M_2$.
As in $\mathcal M_1$, (A1) and (A2) easily hold; here $h_2(\rho_2)=\frac{(\rho_2-\rho_0)^2}{2(1-\rho^2_0)}$
is of the same form as $h_1$. With respect to (A3) we verify pontwise convergence as required, rather than
uniform convergence as in $\mathcal M_1$. Using (\ref{eq:ergodic1}), (\ref{eq:decomp}), (\ref{eq:ergodic2}) 
and (\ref{eq:ergodic3}), it is easily seen that $\frac{\log R_T(\rho_2)}{T}+h_2(\rho_2)\rightarrow 0$ almost surely, 
for all $\rho_2\in\Theta_2$. As in $\mathcal M_1$, it is clear that $\pi(I|\mathcal M_2)=0$ so that (A4) holds.

As regards (A5), note that 
\begin{equation}
h_2\left(\Theta_2\right)=\min\left\{\frac{(1-\rho_0)^2}{2(1-\rho^2_0)},\frac{(1+\rho_0)^2}{2(1-\rho^2_0)}\right\}.
\label{eq:h2_ar1}
\end{equation}
%$\underset{\rho_1\in \Theta_1}{\mbox{ess~inf}}~h_1(\rho_1)
%=\underset{\rho_1\in \Theta_1}{\mbox{ess~inf}}~\frac{(\rho_1-\rho_0)^2}{2(1-\rho^2_0)}=0$.
Now, in contrast with $\mathcal M_1$, here let $\mathcal G_T=\left\{\rho_2\in\Theta_2:|\rho_2|\leq\exp(\beta T)\right\}$, where
$\beta>h_2(\Theta_2)$, with $h_2(\Theta_2)$ given by (\ref{eq:h2_ar1}).
It is easily seen that $\mathcal G_T\rightarrow\Theta_2$ and $h_2(\mathcal G_T)\rightarrow h_2(\Theta_2)$,
as $T\rightarrow\infty$, so that (A5) (1) holds, and (A5) (2) is satisfied by Markov's inequality.
%By Markov's inequality, $\pi(\mathcal G_T)>1-\left(E|\rho_1|\right)\exp\left(-\beta T\right)$.
%That is, (A5) (2) is satisfied. 
Since $\mathcal G_T$ and $\mathbb S$ are compact, verification of (A5) (3) follows in the same way as our proof of 
uniform convergence of $\frac{\log R_T(\cdot)}{T}$
to $-h_1(\cdot)$ in the case of $\mathcal M_1$, provided in Section \ref{subsec:A3}.
That is, (A5) is satisfied for $\mathcal M_2$.

To verify (A6), first note that due to compactness of $\mathcal G_T$, the mean value theorem for integrals yields
\begin{equation}
\frac{1}{m}\log\int_{\mathcal G_T}R_m(\rho_2)\pi(\rho_2|\mathcal M_2)d\rho_2=
\frac{1}{m}\log\left[R_m(\hat\rho_T)\right]+\frac{1}{m}\log\left[\pi(\mathcal G_T|\mathcal M_2)\right],
\label{eq:mvt2}
\end{equation}
for some $\hat\rho_T\in\mathcal G_T$. 
%Compactness of $\mathcal G_T$ also entails that
%\begin{equation}
%h_2(\mathcal G_T)=h_2(\tilde\rho_T),
%\label{eq:h2_GT}
%\end{equation}
%for some $\tilde\rho_T\in\mathcal G_T$.
%\frac{1}{m}\log\int_{|\rho_1|\leq 1}R_m(\rho_1)\pi(\rho_1|\mathcal M_1)d\rho_1$, by the mean value theorem for integrals,
%\begin{equation}
%\frac{1}{m}\log\int_{\mathcal G_T}R_m(\rho_1)\pi(\rho_1|\mathcal M_1)d\rho_1
%=\frac{1}{m}\log \left[R_m(\hat\rho_T)\pi(\Theta_1|\mathcal M_1)\right]
%=\frac{1}{m}\log\left[R_m(\hat\rho_T)\right],
%\label{eq:mvt}
%\end{equation}
%for $\hat\rho_T\in[-1,1]$ depending upon the data. 

%Next observe that if $(-1,1)\subset[-\exp(\beta T),\exp(\beta T)]$, then
%$h_1(\mathcal G_T)=h_1((-1,1))=0$, otherwise, $(-1,1)\supset[-\exp(\beta T),\exp(\beta T)]$, so that
%$\mathcal G_T=[-\exp(\beta T),\exp(\beta T)]$ is compact, and 
%%by compactness of $\tilde{\mathcal G_T}$, 
%hence, $h_1(\mathcal G_T)=h_1(\hat\rho_T)$, for
%$\hat\rho_T\in\mathcal G_T$. 
%That is, 
%\begin{equation}
%h_1(\mathcal G_T)=\max\left\{0,h_1(\hat\rho_T)\right\}. 
%\label{eq:h1_GT}
%\end{equation}
%Note that $h_1(\hat\rho_T)$ may also be zero.

%Since $h_1(\tilde\rho_T)\geq \max\left\{0,h_1(\hat\rho_T)\right\}=h_1(\mathcal G_T)$, 
%$$\frac{1}{m}\log\int_{\tilde{\mathcal G_T}}R_m(\rho_1)\pi(\rho_1|\mathcal M_1)d\rho_1>\delta-h_1(\mathcal G_T)$$ 
%implies, in conjunction with (\ref{eq:mvt}) and (\ref{eq:h1_GT}) that

Since $h_2(\tilde\rho_T)\geq h_2(\mathcal G_T)$, it follows from 
$$\frac{1}{m}\log\int_{\mathcal G_T}R_m(\rho_2)\pi(\rho_2|\mathcal M_2)d\rho_2>\delta-h_2(\mathcal G_T)$$ 
and (\ref{eq:mvt2}) that
$$
\frac{1}{m}\log R_m(\hat\rho_T)+h_2(\hat\rho_T)>\delta-\frac{1}{m}\log\pi(\mathcal G_T|\mathcal M_2)
+h_2(\hat\rho_T)-h_2(\mathcal G_T)>\delta.
$$
The rest of the validation of condition (A6) follows in the same way 
as in the case of $\mathcal M_1$, as detailed in Section \ref{subsec:A6}.

Hence, Theorem 2 %\ref{theorem:bf_convergence} 
of our main manuscript holds, so that
\begin{equation}
\underset{T\rightarrow\infty}{\lim}~\frac{1}{T}\log [B^{(12)}_T]=h_2(\Theta_2),
\label{eq:bf2}
\end{equation}
that is, the Bayes factor heavily favours the (asymptotically) stationary model $\mathcal M_1$ over the nonstationary model $\mathcal M_2$.
Since the true model $P$ is (asymptotically) stationary, this result is very encouraging.

\subsection{Convergence of Bayes factor when $\rho_1$, $\rho_2$, $\sigma_1$ and $\sigma_2$ are all unknown}
\label{subsec:sigma_unknown}

When apart from unknown $\rho_1$ and $\rho_2$, the error variances $\sigma^2_1$ and $\sigma^2_2$ 
associated with models $\mathcal M_1$ and $\mathcal M_2$ are also unknown, 
we consider the parameter spaces $\Theta_1=\left\{(\rho_1,\sigma_1):|\rho_1|<1,\sigma_1\geq \eta\right\}$
and $\Theta_2=\left\{(\rho_2,\sigma_2):\rho_2\in(-1,1)^c\cap\mathbb S,\sigma_2\geq \eta\right\}$ associated with models $\mathcal M_1$
and $\mathcal M_2$, respectively, where $0<\eta<\sigma_0$ is some small constant. For $i=1,2,$ we assume joint priors
$\pi(\rho_i,\sigma_i|\mathcal M_i)$, having densities on $\Theta_i$, with respect to the Lebesgue measure.
It can be easily seen that in this case, for $i=1,2$,
\begin{equation}
h_i(\rho_i,\sigma_i)=\frac{1}{2(1-\rho^2_0)}\left[\left(\rho_0-\frac{\sigma_0\rho_i}{\sigma_i}\right)^2
+\frac{\sigma^2_0}{\sigma^2_i}-(1-\rho^2_0)\log\frac{\sigma^2_0}{\sigma^2_i}-1\right].
\label{eq:h_i}
\end{equation}
Since $(1-\rho^2_0)\log\frac{\sigma^2_0}{\sigma^2_i}+1\leq\log\frac{\sigma^2_0}{\sigma^2_i}+1\leq\frac{\sigma^2_0}{\sigma^2_i}$,
(\ref{eq:h_i}) is non-negative. Also, as in the case with $\sigma_1=\sigma_2=\sigma_0$, it holds that
$h_1(\Theta_1)=0$ and 
$h_2\left(\Theta_2\right)=\min\left\{\frac{(1-\rho_0)^2}{2(1-\rho^2_0)},\frac{(1+\rho_0)^2}{2(1-\rho^2_0)}\right\}$.
Further, note that $\pi(I|\mathcal M_i)=0$, for $i=1,2$.
Thus, conditions (A1)--(A4) are easily seen to hold for both the competing models.

We now verify the remaining conditions for the models.
As regards $\mathcal G_T$, here we set 
%$\mathcal G_T=\left\{(\rho_1,\sigma_1):|\rho_1|<1,\sigma_1\geq 0\right\}\cap
%\left\{(\rho_1,\sigma_1):\sqrt{\rho^2_1+\sigma^2_1}\leq \exp(\beta T)\right\}$ 
$$\mathcal G_T=\left\{(\rho_1,\sigma_1):|\rho_1|<1,\eta\leq\sigma_1\leq\exp(\beta T)\right\}$$
for model $\mathcal M_1$ where $\beta>h_1(\Theta_1)=0$, 
and for model $\mathcal M_2$ we set 
%$\mathcal G_T=\left\{(\rho_2,\sigma_2):|\rho_2|\geq 1,\sigma_2\geq 0\right\}\cap
%\left\{(\rho_2,\sigma_2):\sqrt{\rho^2_2+\sigma^2_2}\leq \exp(\beta T)\right\}$.
%$\mathcal G_T=\left\{(\rho_2,\sigma_2):|\rho_2|\geq 1,\sigma_2\geq 0\right\}\cap
%\left\{(\rho_2,\sigma_2):0\leq\sigma_2\leq \exp(\beta T)\right\}$.
$$\mathcal G_T=\left\{(\rho_2\in\Theta_2,\sigma_2\geq 0):|\rho_2|\leq\exp(\beta T),\eta\leq\sigma_2\leq \exp(\beta T)\right\},$$
where $\beta>h_2(\Theta_2)$.
Note that there exists $T_0\geq 1$ such that $\sigma_0\leq \exp(\beta T)$ for $T\geq T_0$. Hence,
$h_1(\mathcal G_T)=h_1(\Theta_1)=0$ and $h_2(\mathcal G_T)=h_2(\Theta_2)=
\min\left\{\frac{(1-\rho_0)^2}{2(1-\rho^2_0)},\frac{(1+\rho_0)^2}{2(1-\rho^2_0)}\right\}$, for $T\geq T_0$. 
Hence, (A5) (1) holds for both $\mathcal M_1$ and $\mathcal M_2$. 
Now observe that $$\pi(\mathcal G_T|\mathcal M_1)=\pi\left(\eta\leq\sigma_1\leq\exp(\beta T)\right)
>1-E\left(\sigma_1\right)\exp\left(-\beta T\right),$$ so that (A5) (2) holds for $\mathcal M_1$.
For $\mathcal M_2$, denoting by $E_2$ the expectation with respect to $\pi(\cdot|\mathcal M_2)$, note that 
$$\pi(\mathcal G_T|\mathcal M_2)=\pi\left(\eta\leq\sigma_2\leq\exp\left(\beta T\right)|\mathcal M_2\right)
-\pi\left(|\rho_2|>\exp(\beta T),\rho_2\in\mathbb S,\eta\leq\sigma_2\leq\exp(\beta T)|\mathcal M_2\right),$$ where
$$\pi\left(\eta\leq\sigma_2\leq\exp\left(\beta T\right)|\mathcal M_2\right)>1-E_2\left(\sigma_2\right)\exp(-\beta T)$$
and $$\pi\left(|\rho_2|>\exp(\beta T),\rho_2\in\mathbb S,\eta\leq\sigma_2\leq\exp(\beta T)|\mathcal M_2\right)\leq 
\pi\left(|\rho_2|>\exp(\beta T)|\mathcal M_2\right)
<E_2|\rho_2|\exp(-\beta T),$$ by Markov's inequality. It follows that
$$\pi(\mathcal G_T|\mathcal M_2)>1-\left(E_2(\sigma_2)+E_2|\rho_2|\right)\exp(-\beta T),$$ that is, (A5) (2) holds
for $\mathcal M_2$.
That (A5) (3) holds for $\mathcal M_1$
can be shown in the same way as in Section \ref{subsec:A3}, by replacing $|\rho_1|<1$ by $|\rho_1|\leq 1$ in $\mathcal G_T$ and using the assumption that $\sigma_1\geq\eta>0$.
For $\mathcal M_2$ as well, (A5) (3) can be seen to hold in the same way using compactness of $\mathcal G_T$ and $\mathbb S$, and the assumption that $\sigma_2$ is bounded away from zero.

Now observe that for model $\mathcal M_1$, since $h_1(\mathcal G_T)=0$ for $T\geq T_0$, it can be shown in the same way
as in Section \ref{subsec:A6} that
$$\frac{1}{m}\log R_m(\hat\rho_T,\hat\sigma_T)+h_1(\hat\rho_T,\hat\sigma_T)>\delta$$
holds for $T\geq T_0$.
The same holds for model $\mathcal M_2$ using compactness of $\mathcal G_T$, as shown in Section \ref{subsec:verify_M2}
in the context of verification of (A6) for $\mathcal M_2$ when $\sigma_2=\sigma_0$.
Finally observe that it is sufficient to establish convergence of $\sum_{T=T_0}^{\infty}P\left(\tau(\mathcal G_T,\delta)>T\right)$
for large enough $T_0$, which can be done similarly as before, for both $\mathcal M_1$ and $\mathcal M_2$.

%Now note that the continuous function $\psi(\theta)=\|\theta\|=\sqrt{\theta^2_1+\theta^2_2}$ of $\theta=(\theta_1,\theta_2)$ 
%is coercive, that is, for every sequence $\{\theta_T:T>0\}$ such that $\|\theta_T\|\rightarrow\infty$,
%$\psi(\theta_T)\rightarrow\infty$, from which it follows that 
%$\left\{(\rho_1,\sigma_1):\sqrt{\rho^2_1+\sigma^2_1}\leq \exp(\beta T)\right\}$
%is coercive, and hence, compact, for any $T>0$. For model $\mathcal M_1$, if 
%$\left\{(\rho_1,\sigma_1):|\rho_1|<1,\sigma_1\geq 0\right\}\subset 
%\left\{(\rho_2,\sigma_2):\sqrt{\rho^2_2+\sigma^2_2}\leq \exp(\beta T)\right\}$, 
%then $h_1(\mathcal G_T)=0$, otherwise, $\mathcal G_T$ is compact, so that $h_1(\mathcal G_T)
%=\max\left\{0,h_1(\hat\rho_T)\right\}$ holds for some $\hat\rho_T\in\mathcal G_T$,
%as in (\ref{eq:h1_GT}). As regards, model $\mathcal M_2$, observe that the corresponding $\mathcal G_T$ is compact
%for every $T>0$.
%$\left\{\theta\in\mathbb R^2:\psi(\theta)\leq c\right\}$, are compact, for any $c>0$. 
%Hence, the sets $\mathcal G_T$ are intersections with compact sets 
%(for $\mathcal M_1$, we again abuse notation by referring to 
%$\left\{(\rho_1,\sigma_1):|\rho_1|\leq 1,\sigma_1\geq 0,\sqrt{\rho^2_1+\sigma^2_1}\leq \exp(\beta T)\right\}$ as $\mathcal G_T$)
%for every $T>0$. 
%With these, the conditions (A1)--(A6) can be verified as before.

Hence, Theorem 2 of our main manuscript is applicable to this situation and the result remains the same as (\ref{eq:bf2}).

\section{A first look at the applicability of our Bayes factor result to some infinite-dimensional models}
\label{sec:infinite_dimensonal}
\subsection{Traditional Dirichlet process model: undominated case}
\label{subsec:dp_undominated}
Theorem 2 %\ref{theorem:bf_convergence} 
requires the unnormalized posterior to admit factorization 
as the prior times the likelihood. It is well-known that for the original nonparametric models associated with the Dirichlet
process prior (\ctn{Ferguson73}) such factorization is not possible, since there is no parametric form of the
likelihood. In other words, if $[X_1,\ldots,X_T|F]\stackrel{iid}{\sim}F$, where $F\sim DP\left(\alpha F_0\right)$,
where $DP\left(\alpha F_0\right)$ stands for Dirichlet process with base measure $F_0$ and precision parameter $\alpha$,
then the likelihood associated with the data $X_1,\ldots,X_T$ does not have a parametric form, and although the
posterior $\pi(F|\bX_T)$ is well-defined, it is not dominated by any $\sigma$-finite measure
(see, for example, Proposition 7.7 of \ctn{Orbanz14}), and hence does not have a density. This of course prevents
factorization of the posterior of $F$ as the prior times likelihood. Moreover, recall that \ctn{Shalizi09} also assumes
the existence of a common reference measure for the posteriors $\pi(\cdot|\bX_T)$, for all $T$, which does not hold
here. Indeed, such an assumption is valid in the usual dominated case of Bayes theorem where the aforementioned factorization
is possible; in such (usually parametric) cases, the prior is the natural common dominating measure 
(see \ctn{Schervish95}, for example).

\subsection{Dirichlet process mixture model: dominated case}
\label{subsec:dpm_dominated}
Since Dirichlet process supports discrete distributions with probability one, 
the modeling style described in Section \ref{subsec:dp_undominated} is inappropriate if the data 
$\bX_T$ arises from some continuous distribution.
Hence, for such data it is usual in Bayesian nonparametrics based on the Dirichlet process prior to consider the    
following mixture model (see, for example, \ctn{Ghosh03}): 
\begin{equation}
[X_1,\ldots,X_T|F]\stackrel{iid}{\sim}\int f(\cdot|\xi)dF(\xi), 
\label{eq:dp2}
\end{equation}
where
$f(\cdot|\xi)$ is some standard continuous density, usually Gaussian, given $\xi\sim F$, 
where $F\sim DP\left(\alpha F_0\right)$.
By Sethuraman's construction (\ctn{Sethuraman94}), $F(\cdot)=\sum_{i=1}^{\infty}p_i\delta_{\xi_i}(\cdot)$, with probability
one, where, for $i=1,2,\ldots$, $\xi_i\stackrel{iid}{\sim}F_0$, and for any $\xi$, $\delta_{\xi}(\cdot)$ 
denotes the point mass on $\xi$. Also, for $i=1,2,\ldots$,
$p_i=V_i\prod_{j<i}(1-V_j)$, where $V_i\stackrel{iid}{\sim}Beta(\alpha,1)$.
It is easy to verify that $\sum_{i=1}^{\infty}p_i=1$, almost surely.
Application of Sethuraman's construction in (\ref{eq:dp2}) yields the equivalent infinite mixture representation
\begin{equation}
[X_1,\ldots,X_T|\theta]\stackrel{iid}{\sim}\sum_{i=1}^{\infty}p_if(\cdot|\xi_i), 
\label{eq:dp3}
\end{equation}
where $\theta=(\xi_1,\xi_2,\ldots,V_1,V_2,\ldots)$ is the infinite-dimensional parameter. The prior
on $\theta$ is already specified by the $iid$ $F_0$ and $Beta(\alpha,1)$ distributions, and is the infinite product
probability measure associated with these $iid$ distributions, so that each factor of the product of the
probability measures is dominated by the Lebesgue measure.
In this case, the posterior of $\theta$ admits the representation
\begin{equation}
\pi(\theta|\bX_T)\propto \pi(\theta)\prod_{t=1}^T\left[\sum_{i=1}^{\infty}p_if(X_t|\xi_i)\right],
\label{eq:dp4}
\end{equation}
and hence the representation of Bayes factor in terms of the prior and the likelihood holds in this case, as required
by Theorem 2 of our main manuscript. %\ref{theorem:bf_convergence}. 
Moreover, the posterior $\pi(\cdot|\bX_T)$ is absolutely continuous with 
respect to $\pi(\cdot)$ for all $T$, as assumed by \ctn{Shalizi09}.

\subsection{Polya urn based mixture obtained by integrating out random $F$: dominated case but $\mathcal T$ changes
with $T$}
\label{subsec:polya_undominated}
Assume that for $t=1,\ldots,T$, $[X_t|\phi_t]\sim f(\cdot\phi_t)$, independently, and $\phi_1,\ldots,\phi_T\stackrel{iid}{\sim}F$,
where $F\sim DP\left(\alpha F_0\right)$. This is equivalent to the Dirichlet process mixture model (\ref{eq:dp2}),
but if $F$ is integrated out, then the joint distribution of $\phi_1,\ldots,\phi_T$ is given by the Polya urn scheme,
that is, $\phi_1\sim F_0$, and for $t=2,\ldots,T$, $[\phi_t|\phi_1,\ldots,\phi_{t-1}]\sim\frac{\alpha F_0}{\alpha+t-1}
+\frac{\sum_{j=1}^{t-1}\delta_{\phi_j}}{\alpha+t-1}$ (see, for example, \ctn{Ferguson73}, \ctn{Escobar95}). 
The joint prior distribution of $\phi_1,\ldots,\phi_T$
has a density with respect to a measure composed of Lebesgue measures in lower dimensions; see Lemma 1.99 of \ctn{Schervish95}
for the exact forms of the density and the dominating measure. Hence, in this case the posterior of 
$\phi_1,\ldots,\phi_T$ is proportional to the prior times the likelihood, where the likelihood is given by 
$\prod_{t=1}^Tf(X_t|\phi_t)$, and the posterior is dominated by the prior probability measure. Hence, a countably infinite convex
combination of the prior probability measures dominates the posterior of $\phi_1,\ldots,\phi_T$ for all $T$, as
required for the results of \ctn{Shalizi09} to hold. However,
\ctn{Shalizi09} assumes that the $\sigma$-field $\mathcal T$ associated with the parameter space $\Theta$ does not change
with $T$, which does not hold in this case.

\subsection{Polya urn based finite mixture: dominated case and $\mathcal T$ remains fixed}
\label{subsec:sb_dominated}

\ctn{Bhattacharya08} (see also \ctn{Sabya11}, \ctn{Sabya12}) introduce the following finite mixture model based 
on Dirichlet process:
\begin{align}
X_1,\ldots,X_T&\stackrel{iid}{\sim}\frac{1}{M}\sum_{i=1}^Mf(\cdot|\phi_i);\label{eq:sb1}\\
\phi_1,\ldots,\phi_M&\stackrel{iid}{\sim}F;\label{eq:sb2}\\
F&\sim DP\left(\alpha F_0\right),\label{eq:sb3}
\end{align}
where $f(\cdot|\phi)$ is any standard density as before, given parameter(s) $\phi$, and $M~(>1)$ is some fixed integer.
Integrating out $F$ yields the following Polya urn scheme for the joint distribution of $\phi_1,\ldots,\phi_M$:
$\phi_1\sim F_0$, and for $t=2,\ldots,M$, $[\phi_t|\phi_1,\ldots,\phi_{t-1}]\sim\frac{\alpha F_0}{\alpha+t-1}
+\frac{\sum_{j=1}^{t-1}\delta_{\phi_j}}{\alpha+t-1}$. Here $\theta=(\phi_1,\ldots,\phi_M)$, which is of fixed, finite size,
even though the problem is induced by the nonparametric Dirichlet process prior. Also clearly the
$\sigma$-field $\mathcal T$ associated with the parameter space $\Theta$ does not change with $T$. 
Thus, in this set-up, not only is the posterior
written in terms of product of the prior and the likelihood, but is dominated by the Polya urn based prior of $\theta$,
for all sample sizes $T$.

\subsection{Nonparametric Bayesian using the Polya tree prior: dominated case}
\label{subsec:polya_tree_dominated}
\ctn{Lavine92}, \ctn{Lavine94} proposed the Polya tree prior for the random probability measure $F$ 
as an alternative to the Dirichlet process prior.
Briefly, one starts with a partition $\pi_1=\{B_0,B_1\}$ of the sample space $\Omega$, so that $\Omega=B_0\cup B_1$.
This procedure is then continued with $B_0=B_{00}\cup B_{01}$, $B_1=B_{10}\cup B_{11}$, etc. At level $m$, the
partition is then $\pi_m=\{B_{\epsilon}:\epsilon=\epsilon_1\ldots\epsilon_m\}$, where $\epsilon$ are all binary
sequences of length $m$. Let $\Pi=\left\{\pi_m:m=1,2,\ldots\right\}$, and $\mathcal A=\left\{\alpha_{\epsilon}\right\}$
be a sequence of non-negative numbers, one for each partitioning subset.
Now, if $Y_{\epsilon 0}=F(B_{\epsilon 0}|B_{\epsilon})\sim Beta(\alpha_{\epsilon 0},\alpha_{\epsilon 1})$
independently with respect to the $\epsilon$'s, then $F$
is said to have the Polya tree prior $PT(\Pi,\mathcal A)$. 

It can be shown that if $\alpha_{\epsilon}\propto m^{-1/2}$, the Polya tree prior reduces to the Dirichlet process prior,
confirming that the latter is a special case of the Polya tree prior. However, the most important property of the 
Polya tree prior is that with appropriate choices of the $\alpha_{\epsilon}$, 
$F$ can be made absolutely continuous with respect to the Lebesgue measure. Specifically, if $\alpha_{\epsilon}\propto m^2$, 
for the $m$-th level subset, then $F$ is dominated by the Lebesgue measure almost surely. 
Hence, if $[X_1,\ldots,X_T|F]\sim F$ and $F\sim PT(\Pi,\mathcal A)$, with $\alpha_{\epsilon}\propto m^2$, then 
the likelihood is available almost surely. Here we may set 
$\theta=\left\{Y_{\epsilon 0}:\epsilon=\epsilon_1\ldots\epsilon_m,m=1,2,\ldots\right\}$, which has the infinite product
prior measure. The posterior of $F$ given $\bX_T$, which is also a Polya tree process, is dominated by $\pi(\theta)$
for all $T>0$. Similar issues hold for the extended Polya tree prior, namely, the optional Polya tree prior proposed
by \ctn{Wong10}.

\subsection{Bayesian density estimation using the generalized lognormal process prior: dominated case}
\label{subsec:lognormal_dominated}
\ctn{Lenk88} model the unknown density $f(x)$ with respect to measure $\lambda$ as
\begin{equation}
f(x)=\frac{W(x)}{\int_{\mathcal X}W(s)d\lambda(s)},
\label{eq:lognormal1}
\end{equation}
where $W$ is a generalized lognormal process over $\mathcal X$. The generalized lognormal process has distribution
$\Lambda_{\eta}$ given by (see \ctn{Lenk88})
\begin{align}
\Lambda_{\eta}(A)&=\frac{E\left[\left(\int_{\mathcal X}Wd\lambda\right)^{\eta}\mathbb I_A\right]}
{E\left[\left(\int_{\mathcal X}Wd\lambda\right)^{\eta}\right]},
\label{eq:lognormal2}
\end{align}
where $-\infty<\eta<\infty$ and the expectations are taken with respect to the usual lognormal process, that is, with respect to 
$W=\exp\left(Z\right)$, where $Z$ is a Gaussian process. In (\ref{eq:lognormal2}), $\mathbb I_A$ is the indicator
of the set $A$, where $A$ belongs to the Borel $\sigma$-field associated with the space of functions from $\mathcal X$ to
$(0,\infty)$. The properties and moments of the lognormal process are provided in \ctn{Lenk88}. 

In this formulation, the likelihood with respect to $iid$ data $X_1,\ldots,X_T$ is defined via (\ref{eq:lognormal1}). 
The prior distribution, as well as the posterior distribution of $\Theta=W$ for all $T\geq 1$, 
are absolutely continuous with respect to the distribution of the lognormal process $W=\exp\left(Z\right)$, 
where $Z$ is a Gaussian process.

\subsection{Bayesian regression using Gaussian process: dominated case}
\label{subsec:gp_dominated}
Consider the following regression model with covariates $\left\{C_t:t=1,\ldots,T\right\}$ (see \ctn{Choi07}, for example):
\begin{align}
X_t&=\zeta(C_t)+\epsilon_t,~t=1,\ldots,T;\notag\\
\epsilon_t&\stackrel{iid}{\sim}N\left(0,\sigma^2\right);\notag\\
\sigma&\sim \varphi;\notag\\
\zeta(\cdot) &\sim GP\left(\mu(\cdot),K(\cdot,\cdot)\right),
\label{eq:gp1}
\end{align}
where $\varphi$ is a probability measure on the positive part of the real line, and in (\ref{eq:gp1}), 
$GP\left(\mu(\cdot),K(\cdot,\cdot)\right)$ stands for the Gaussian process with mean function $E\left[\zeta(c)\right]=\mu(c)$ 
for all $c\in\mathfrak C$, where $\mathfrak C$ is the space of covariates, and positive definite covariance function
$Cov\left(\zeta(c_1),\zeta(c_2)\right)=K(c_1,c_2)$, for all $c_1,c_2\in\mathcal C$.
Here, by positive definite function $K(\cdot,\cdot)$ on $\mathfrak C\times\mathfrak C$, we mean
$\int K(c,c')g(c)g(c')d\nu(c)d\nu(c')>0$ for all non-zero functions $g\in L_2\left(\mathfrak C,\nu\right)$,
where $L_2\left(\mathfrak C,\nu\right)$ denotes the space of functions square-integrable on $\mathfrak C$ with respect to
the measure $\nu$.

In what follows we borrow the statements of the following definition of eigenvalue and eigenfunction, and the subsequent
statement of Mercer's theorem from \ctn{Rasmussen06}.
\begin{definition}
\label{def:eigen}
A function $\psi(\cdot)$ that obeys the integral equation
\begin{equation}
\int_{\mathfrak C}K(c,c')\psi(c) d\nu(c) = \lambda\psi(c'),
\label{eq:eigen}
\end{equation}
is called an eigenfunction of the kernel $K$ with eigenvalue $\lambda$ with respect to the 
measure $\nu$. 
\end{definition}

We assume that the ordering is chosen such that $\lambda_1\geq\lambda_2\geq\cdots$.
The eigenfunctions are orthogonal with respect to $\nu$ and can be chosen to be
normalized so that $\int_{\mathfrak C}\psi_i(c)\psi_j(\bx)d\nu(c)=\delta_{ij}$,
where $\delta_{ij}=1$ if $i=j$ and $0$ otherwise.

The following well-known theorem (see, for example, \ctn{Konig86}) expresses the positive definite kernel 
$K$ in terms of its eigenvalues and eigenfunctions.

\begin{theorem}[Mercer's theorem]
\label{theorem:mercer}
Let $(\mathfrak C,\nu)$ be a finite measure space and $\mathcal C\in L_{\infty}\left(\mathfrak C^2,\nu^2\right)$
be a positive definite kernel. By $L_{\infty}\left(\mathfrak C^2,\nu^2\right)$ we mean
the set of all measurable functions $K:\mathfrak C^2\mapsto\mathbb R$ which are essentially bounded,
that is, bounded up to a set of $\nu^2$-measure zero. For any function $K$ in this set, its essential supremum,
given by $\inf\left\{r\geq 0:|\mathcal K(c,c')|<r,~\mbox{for almost all}~(c,c')\in\mathfrak C\times\mathfrak C\right\}$
serves as the norm $\|K\|$.
%such that $T_K : L_2\left(\mathfrak X,\mu\right)\mapsto L_2\left(\mathfrak X,\mu\right)$ is positive
%definite. 

Let $\psi_j\in L_2\left(\mathfrak C,\nu\right)$ be the normalized eigenfunctions of
$K$ associated with the eigenvalues $\lambda_j(K)>0$. Then
\begin{itemize}
\item[(a)] the eigenvalues $\left\{\lambda_j(K)\right\}_{j=1}^{\infty}$ are absolutely summable.
\item[(b)] $K(c,c')=\sum_{j=1}^{\infty}\lambda_j(K)\psi_j(c)\bar{\psi_j}(c')$
holds $\nu^2$-almost everywhere, where the series converges absolutely and uniformly
$\nu^2$-almost everywhere. In the above,
$\bar{\psi_j}$ denotes the complex conjugate of $\psi_j$.
\end{itemize}
\end{theorem}

It follows that the Gaussian process $\zeta$ admits the representation below almost surely:
\begin{equation}
\zeta(\cdot)=\mu(\cdot)+\sum_{i=1}^{\infty}\sqrt{\lambda_i}\psi_i(\cdot)e_i,
\label{eq:kl1}
\end{equation}
where, for $i=1,2,\ldots$, $e_i\stackrel{iid}{\sim}N(0,1)$. The above representation for Gaussian processes
is popularly known as the Karhunen-Lo\`{e}ve expansion (see, for example, \ctn{Ash75}).

Hence, both the likelihood and the prior can be parameterized in terms of $\psi_i(\cdot);~i=1,2,\ldots$
and $\be=\left\{e_i;~i=1,2,\ldots\right\}$, the latter being unknown and having the infinite product prior distribution 
such that $e_i\stackrel{iid}{\sim}N(0,1)$; $i=1,2,\ldots$. Letting $\theta=(\be,\sigma)$, note that 
the posterior $\pi(\theta|\bX_T)$, for all $T>0$, 
is clearly dominated by this infinite product prior measure times $\varphi$.

\renewcommand\baselinestretch{2}
\normalsize
\bibliographystyle{ECA_jasa}
\bibliography{irmcmc}

\end{document}